\documentclass[12pt]{amsart}
\usepackage{amsmath,amscd,amssymb,amsfonts,graphics}
\usepackage{xcolor,hyperref}
\hypersetup{colorlinks=true,citecolor=black,linkcolor=black}
\newcommand{\hl}{\hyperlink}
\newcommand{\htt}{\hypertarget}
\newcommand{\h}{\hbox}
\newcommand{\q}{\quad}
\newcommand{\nin}{\noindent}
\newcommand{\bs}{\par\bigskip}
\newcommand{\ms}{\par\medskip}
\newcommand{\sk}{\par\smallskip}
\newcommand{\bsn}{\par\bigskip\noindent}
\newcommand{\msn}{\par\medskip\noindent}
\newcommand{\skn}{\par\smallskip\noindent}
\newcommand{\ges}{\geqslant}
\newcommand{\les}{\leqslant}
\newcommand{\1}{\hskip1pt}
\newcommand{\mcap}{\hbox{$\bigcap$}}
\newcommand{\mcup}{\hbox{$\bigcup$}}
\newcommand{\msum}{\hbox{$\sum$}}
\newcommand{\mopl}{\hbox{$\bigoplus$}}
\newcommand{\mprod}{\hbox{$\prod$}}

\newcommand{\K}{{\mathcal K}}

\newcommand{\M}{{\mathcal M}}
\newcommand{\Nc}{{\mathcal N}}
\newcommand{\OO}{{\mathcal O}}

\newcommand{\Q}{{\mathbb Q}}
\newcommand{\C}{{\mathbb C}}
\newcommand{\N}{{\mathbb N}}
\newcommand{\R}{{\mathbb R}}

\newcommand{\Z}{{\mathbb Z}}
\newcommand{\ob}{{\mathbf 1}}
\newcommand{\Gr}{{\rm Gr}}
\newcommand{\al}{\alpha}
\newcommand{\be}{\beta}
\newcommand{\ga}{\gamma}
\newcommand{\Ga}{\Gamma}
\newcommand{\de}{\delta}
\newcommand{\De}{\Delta}
\newcommand{\la}{\lambda}
\newcommand{\si}{\sigma}

\newcommand{\ep}{\varepsilon}
\newcommand{\qh}{\widehat{q}}

\newcommand{\nut}{\widetilde{\nu}}
\newcommand{\gat}{\widetilde{\gamma}}
\newcommand{\Ht}{{}\,\widetilde{\!H}{}}

\newcommand{\Om}{\Omega}
\newcommand{\dd}{\partial}
\newcommand{\ddd}{{\rm d}}
\newcommand{\ee}{{\bf e}}
\newcommand{\ff}{{\bf f}}
\newcommand{\Ff}{F_{\!f}}
\newcommand{\Gp}{\Gamma_{\!+}}
\newcommand{\Gf}{\Gamma_{\!f}}
\newcommand{\Gg}{\Gamma_{\!g}}
\newcommand{\Spx}{{\rm Supp}_{(x)}}
\newcommand{\CFf}{{\rm CF}_{\hskip-2pt f}}
\newcommand{\CFfb}{{\rm CF}_{\hskip-2pt f,{\bf b}}}
\newcommand{\CFfbi}{{\rm CF}_{\hskip-2pt f,{\bf b},i}}
\newcommand{\CFfin}{{\rm CF}_{\hskip-2pt f,{\rm in}}}
\newcommand{\CFfla}{{\rm CF}_{\hskip-2pt f,\la}}
\newcommand{\CFfinla}{{\rm CF}_{\hskip-2pt f,{\rm in},\la}}
\newcommand{\CFfl}{{\rm CF}_{\hskip-2pt f+\ell^r}}
\newcommand{\CFflc}{{\rm CF}_{\hskip-2pt f+\ell^r\!,{\bf c}}}
\newcommand{\CFflci}{{\rm CF}_{\hskip-2pt f+\ell^r\!,{\bf c},i}}
\newcommand{\CFflinc}{{\rm CF}_{\hskip-2pt f+\ell^r\!,{\rm in},{\bf c}}}
\newcommand{\CFflinci}{{\rm CF}_{\hskip-2pt f+\ell^r\!,{\rm in},{\bf c},i}}
\newcommand{\Sp}{{\rm Sp}}
\newcommand{\SpW}{{}^W\!\1{\rm Sp}}
\newcommand{\SpG}{{\rm Sp}^{\Gamma}}
\newcommand{\SpGp}{{\rm Sp}^{\Gamma\1\prime}}
\newcommand{\SpDef}{{\rm Sp}^{\rm Def}}
\newcommand{\SpWG}{{}^W{\rm Sp}^{\Gamma}}
\newcommand{\SpDefW}{{}^W{\rm Sp}^{\rm Def}}

\newcommand{\tos}{\,{\to}\,}
\newcommand{\eq}{\,{=}\,}
\newcommand{\defs}{\,{:=}\,}
\newcommand{\nes}{\,{\ne}\,}
\newcommand{\ins}{\,{\in}\,}
\newcommand{\sst}{\,{\subset}\,}
\newcommand{\stm}{\,{\setminus}\,}
\newcommand{\gess}{\,{\ges}\,}
\newcommand{\less}{\,{\les}\,}
\newcommand{\sgt}{\,{>}\,}
\newcommand{\slt}{\,{<}\,}
\newcommand{\pl}{\1{+}\1}
\newcommand{\mi}{\1{-}\1}
\newcommand{\bl}{\bigl}
\newcommand{\br}{\bigr}

\newcommand{\ssb}{\raise.15ex\h{${\scriptscriptstyle\bullet}$}}
\newcommand{\ssc}{\,\raise.15ex\h{${\scriptstyle\circ}$}\,}

\newcommand{\into}{\hookrightarrow}
\newcommand{\simto}{\,\,\rlap{\hskip1.5mm\raise1.4mm\hbox{$\sim$}}\hbox{$\longrightarrow$}\,\,}

\begin{document}
\title[Spectrum of non-degenerate functions]
{Spectrum of non-degenerate functions\\with simplicial Newton polytopes}
\author[S.-J. Jung]{Seung-Jo Jung}
\address{S.-J. Jung : Department of Mathematics Education, and Institute of Pure and Applied Mathematics, Jeonbuk National University, Jeonju, 54896, Korea}
\email{seungjo@jbnu.ac.kr}
\author[I.-K. Kim]{In-Kyun Kim}
\address{I.-K. Kim : Department of Mathematics, Yonsei University, Seoul, Korea}
\email{soulcraw@gmail.com}
\author[M. Saito]{Morihiko Saito}
\address{M. Saito : RIMS Kyoto University, Kyoto 606-8502 Japan}
\email{msaito@kurims.kyoto-u.ac.jp}
\author[Y. Yoon]{Youngho Yoon}
\address{Y. Yoon : Department of Mathematics, Chungbuk National University, Cheongju-si, Chungcheongbuk-do, 28644, Korea}
\email{mathyyoon@gmail.com}
\thanks{This work was partially supported by NRF grant funded by the Korea government(MSIT) (the first author: NRF-2021R1C1C1004097, the second author: NRF-2020R1A2C4002510, and the fourth author: NRF-2020R1C1C1A01006782). The first author was partially supported by “Research Base Construction Fund Support Program” funded by Jeonbuk National University in 2020. The third author was partially supported by JSPS Kakenhi 15K04816.}
\begin{abstract} We show a precise proof of Steenbrink's formula for the spectrum of convenient Newton non-degenerate functions, and prove the symmetry of combinatorial polynomials in the simplicial case. Combined with the  modified Steenbrink conjecture for spectral pairs (that  is, weighted spectrum) which is recently proved in that case, this simplifies quite a lot of their calculations in such a case. We also introduce the $\Gamma$-spectrum of simplicial convenient non-degenerate functions as a first approximation of the spectrum, generalizing Arnold's picture in the 2 variable case. Analyzing their difference, we can find simple formulas for weighted spectrum in the 3 or 4 variable case. This is proved by using the symmetry of combinatorial polynomials, and fails in the non-simplicial case. Combining these with the Yomdin-Steenbrink formula for the spectrum, we can prove a formula for the spectrum of certain non-isolated surface singularities with simplicial non-degenerate Newton boundaries. As a byproduct of these arguments, we find an example where the Yomdin-Steenbrink formula for spectral pairs does not hold because of fusions of compact faces under projections.
\end{abstract}
\maketitle
\centerline{\bf Introduction}
\bsn
For a holomorphic function $f\,{:}\,(\C^n,0)\tos(\C,0)$, the {\it spectrum\1} $\Sp_f(t)$ is a fractional power polynomial defined from the Hodge filtration with the semi-simple part of the monodromy on the vanishing cohomology, see \cite{St1}, \cite{St2} (and \hl{1.1}{1.1} below). Adding the information of {\it weights\1} in the isolated singularity case, we can define also the {\it weighted spectrum}
$$\SpW_f(t,u)\eq\msum_{i=1}^{\mu_f}\,t^{\al_i}u^{w_i},$$
\par\nin as the generating polynomial of {\it spectral pairs\1} $\{(\al_i,w_i)\}_{i\in[1,\1\mu_f]}$, see \hl{3.2}{3.2} below.
These are quite useful invariants, and have been studied a lot by many people (see for instance \cite{Ar}, \cite{Bu}, \cite{BS}, \cite{DS}, \cite{KL1}, \cite{exp}, \cite{ste}, \cite{Va2}, \cite{Yo} among others). It is, however, not necessarily easy to calculate concretely them for each example even though there is an explicit formula for spectrum in the convenient Newton non-degenerate case, see \cite{St1}, \cite{exp}. Indeed, even in the 3 variable case, it is quite nontrivial to determine the Newton polyhedron (see \hl{2.2}{2.2} below), and moreover we have to calculate the {\it combinatorial polynomials\1} especially for the weighted spectrum as is explained after Theorem\,\,\hl{T1}{1} below.
\sk
In the convenient Newton non-degenerate case with 2 variables, there is a very simple formula for the weighted spectrum, which is proved in \cite{St1}, and is expressed by a picture as in \cite{Ar}, see also Example~\hl{E1.4b}{1.4b} below. One can naturally imagine its generalization to the higher dimensional case assuming that the Newton polytope $\Gp(f)$ is convenient and simplicial. (In this paper, $\Gp(f)$ (or $f$) is called {\it simplicial\1} if any compact face of $\Gp(f)$ is a simplex. This condition does not seem fatally restrictive for explicit calculations of examples. Indeed, if $r$ points are given arbitrarily in $\{1,\dots,m\}^n$ with $r$ sufficiently small compared with $m$, then the possibility of more than $n$ points lying on one hyperplane is quite small.)
\sk
In the simplicial case, one can define the $\Ga$-{\it spectrum\1} $\SpG_f(t)$ from the Newton polytope $\Gp(f)$ using \cite{exp} so that the $\Ga$-spectral numbers have a {\it symmetry with center} $\tfrac{n}{2}$, that is,
\htt{1}{}
$$\SpG_f(t^{-1})\1t^n\eq\SpG_f(t).
\leqno(1)$$
\par\nin Moreover it coincides with the spectrum $\Sp_f(t)$ for the spectral numbers at most 1, that is,
\htt{2}{}
$$\SpG_f(t)^{\les 1}\eq\Sp_f(t)^{\les 1},
\leqno(2)$$
\par\nin with $\Sp_f(t)^{\les 1}\defs\msum_{\al_i\les1}\,t^{\al_i}$ for $\Sp_f(t)\eq\msum_i\,t^{\al_i}$ (where the $\al_i$ are called the {\it spectral numbers}), see \hl{3.1}{3.1} below. In the {\it non-simplicial\1} case, neither the symmetry (\hl{1}{1}) nor the property (\hl{2}{2}) holds even though the $\Ga$-spectrum can be defined as in \cite{St1}. This is closely related to a conjecture of Steenbrink on spectral pairs, see Remark\,\,\hl{R3.1c}{3.1c} below.
\sk
Calculating some examples in the 3 variable case, however, one sees immediately that $\SpG_f(t)\ne\Sp_f(t)$, and the {\it defect\1} $\,\SpDef_f(t)\defs\Sp_f(t)\mi\SpG_f(t)$ is always a sum of fractional power polynomials associated with certain (not necessarily interior) {\it vertices\1} (or 0-dimensional faces) of the Newton polytope $\Gp(f)$. Since the defective spectral numbers are contained in the open interval $(1,2)$ as a consequence of (\hl{1}{1}) and (\hl{2}{2}) when $n\eq3$, the defect can be determined by comparing the $\Ga$-spectrum (modulo the action of multiplication by $t$ on the Laurent fractional polynomials) with the monodromy zeta function \cite{Va1}.
\sk
Let $\CFf^k$ be the set of $k$-dimensional {\it compact faces\1} of $\Gp(f)$. We have the following.
\par\htt{T1}{}\msn
{\bf Theorem\,\,1.} {\it Assume $n\eq 3$, and $f$ is Newton non-degenerate, convenient, and simplicial. Let $\ga_{\si}$ be the number of $2$-dimensional faces $($not necessarily \h{\it compact$\1)$} of $\Gp(f)$ containing $\si\ins\CFf^0$. Then
\htt{3}{}
$$\SpDef_f(t)=\msum_{\si\in\CFf^0}\,(\ga_{\si}{-}3)\1t\1q_{\si}(t),
\leqno(3)$$
\par\nin where $\,q_{\si}(t)\eq\msum_{k=1}^{\de_{\si}-1}\,t\1^{k/\de_{\si}}$ with $\de_{\si}\defs{\rm GCD}(a,b,c)$ if $\,\si\eq\{(a,b,c)\}\sst\N^3$.}
\ms
This immediately implies the following.
\par\htt{C1}{}\msn
{\bf Corollary~1.} {\it With the notation and assumptions of Theorem\,\,$\hl{T1}{1}$, we have $\SpDef_f(t)\eq0$, that is, $\SpG_f(t)\eq\Sp_f(t)$ if and only if $\de_{\si}\eq1$ for any $\si\ins\CFf^0$ with $\ga_{\si}\sgt3$.}
\ms
By the definition of $\Ga$-spectrum (see (\hl{3.1.1}{3.1.1}) below), the equality (\hl{3}{3}) means that
\htt{4}{}
$$\aligned\Sp_f(t)&=\msum_{\si\in\CFfin}\,\bl(\msum_{j=0}^{2-d_{\si}}\,t^j\br)q_{\si}(t)+\bl|\CFfin^0\br|\,\bl(t{+}t^2\br)\\&\q\q{}+{}\msum_{\si\in\CFf^0}\,(\ga_{\si}{-}3)\1t\1q_{\si}(t).\endaligned
\leqno(4)$$
\par\nin Here $\CFfin\eq\bigsqcup_k\CFfin^k$ with $\CFfin^k$ the set of $k$-dimensional {\it interior\1} compact faces $\si$ of $\Gp(f)$ (where interior means that $\si$ is not contained in any coordinate plane of $\R^3$), $d_{\si}\defs\dim\si$, and $q_{\si}(t)$ is the Poincar\'e polynomial of the graded vector space spanned by the monomials $x^{\nu}$ with $\nu\ins\N^3$ contained in the strict interior of the {\it parallelotope\1} $E_{\si}$ spanned by the vertices $v_j$ of the simplex $\si$, that is, $E_{\si}\eq\bl\{\msum_j\,r_jv_j\mid r_j\ins[0,1]\br\}$, where the grading is given so that any vertex $v_j$ of $\si$ has degree 1, see \hl{3.1}{3.1} below for more details.
\sk
Theorem\,\,\hl{T1}{1} easily follows from the symmetry of the combinatorial polynomials (see (\hl{6}{6}) and Theorem\,\,\hl{TA}{A} in Appendix below) in Steenbrink's formula \cite[Theorem\,\,5.7]{St1} (see also (\hl{2.1}{2.1}) below) for the Poincar\'e polynomial of the Newton filtration on the Jacobian ring $B_{\Ga}\defs\Gr_{V_N}^{\ssb}\bl(\C\{x\}/(\dd f)\br)$\,:
\htt{5}{}
$$\aligned p_{B_{\Ga}}(t)&=\msum_{\emptyset\les\si\les\Gf}\,(-1)^{n-d(\si)}\,(1{-}t)^{k(\si)-d(\si)}\1\qh_{\si}(t)\\&=\msum_{\emptyset\les\si\les\Gf}\,r_{\!\si}(t)q_{\si}(t)\q\q\q\q\h{with}\\ \qh_{\si}(t)&:=\msum_{\emptyset\les\tau\les\si}\,q_{\tau}(t),\\ r_{\!\si}(t)&:=\msum_{\si\les\tau\les\Gf}\,(-1)^{n-k(\tau)}\,(t{-}1)^{k(\tau)-d(\tau)},\endaligned
\leqno(5)$$
\par\nin see also \hl{2.1}{2.1} below. Here $\Gf$ denotes the union of compact faces of $\Gp(f)$ with $\,n\mi k(\si)$ the number of coordinate hyperplanes containing $\si$ and $\,d(\si)\defs d_{\si}{+}1$, where $\,k(\emptyset)\eq d(\emptyset)\eq 0$. The summations $\msum_{\emptyset\les\si\les\Gf}$ in (\hl{5}{5}) are taken over any faces of $\Gf$ {\it including\1} $\emptyset$. We call $r_{\!\si}(t)$ the {\it combinatorial polynomial}. The assertion~(\hl{5}{5}) holds {\it without\1} assuming $f$ is simplicial. Note that (\hl{5}{5}) implies a formula for spectrum by \cite{exp}.
\sk
The {\it symmetry of combinatorial polynomials\1} in the simplicial convenient case has been proved quite recently, see Theorem\,\,\hl{TA}{A} in Appendix.
(This theorem is moved to Appendix in this paper from \cite{des}, since it sits better here.) It asserts the {\it symmetries}
\htt{6}{}
$$r_{\!\si}(t)\eq r_{\!\si}(t^{-1})\1t^{n-d(\si)},\,\,\,\,r^{\rm Def}_{\!\si}(t)\eq r^{\rm Def}_{\!\si}(t^{-1})\1t^{n-d(\si)}\,\,\,\,(\emptyset\less\si\less\Gf),
\leqno(6)$$
\par\nin with $r^{\rm Def}_{\!\si}(t)\defs r_{\!\si}(t)\mi r^{\Ga}_{\!\si}(t)$. Here we have by definition
$$r^{\Ga}_{\!\si}(t)\defs\begin{cases}\msum_{j=0}^{n-d(\si)}t^j&(\si\ins\CFfin),\\|\CFfin^0|\,\msum_{j=1}^{n-1}\,t^j&(\si\eq\emptyset),\\ \10&(\h{otherwise}),\end{cases}$$
\par\nin see (\hl{3.1.1}{3.1.1}) below. From this we can deduce the inequalities
\htt{7}{}
$$\deg r^{\rm Def}_{\!\si}(t)\less n\mi d(\si)\mi 1,\q\deg r^{\rm Def}_{\emptyset}(t)\less n{-}2.
\leqno(7)$$
\par\nin \sk
The symmetry (\hl{6}{6}) simplifies very much the computation of $r_{\!\si}(t)$, since the calculation of coefficients of $r_{\!\si}(t)$ becomes easier as $d(\si)$ {\it increases.} Together with (\hl{7}{7}) it implies for instance in the case $n\eq 3$ that we have $r^{\rm Def}_{\!\si}(t)\eq0$ ($d_{\si}\gess 1$) and moreover $r_{\emptyset}(t)$ coincides with the second term of the right-hand side of (\hl{4}{4}), since $q_{\emptyset}(t)\eq 1$.
\sk
It is much simpler to calculate concrete examples using Theorem\,\,\hl{T1}{1} than the Steenbrink formula (\hl{5}{5}), since the former needs only the $\de_{\si}$ for vertices $\si\ins\CFf^0$ with $\ga_{\si}\sgt3$. Note that $\ga_{\si}$ is equal to the number of 1-dimensional faces (not necessarily compact) containing the vertex $P\ins\Gp(f)$ with $\si\eq\{P\}$. On the other hand, the latter requires a computation of $\qh_{\si}(t)$ or $r_{\!\si}(t)$ for {\it any\1} compact faces $\si\less\Gp(f)$ including the case $\si\eq\emptyset$ (although (\hl{6}{6}) and (\hl{7}{7}) simplify certain calculations as is explained above).
Note, however, that the computation of the $q_{\si}(t)$ can be rather non-trivial even in the simplicial case.
\sk
The combinatorial polynomials $r_{\!\si}(t)$ are needed in an essential way for the (modified) {\it Steenbrink conjecture on spectral pairs,} see \hl{3.2}{3.2} and Theorem\,\,\hl{T2}{2} below.
It is rather surprising that all the $r_{\!\si}(t)$ can be controlled only by using the $\ga_{\si}$ with $\dim\si\eq0$, when $n\eq3$. Here we have $\ga_{\si}\sgt3$ if and only if the cone of the dual fan of $\Gp(f)$ corresponding to $\si\ins\CFf^0$ is {\it non-simplicial} (in this case the toric variety can be {\it non-quasi-smooth}).
\sk
For explicit computations of examples, note, however, that it is not necessarily easy to {\it determine\1} the Newton polytope $\Gp(f)$ from the support of $f$ even in the 3 variable case, see (\hl{2.2}{2.2}) below.
This seems almost impossible for $n\eq 4$ without using a computer except for rather simple cases, see for instance \hl{2.4}{2.4} below.
\sk
Theorem\,\,\hl{T1}{1} is closely related to the (modified) {\it Steenbrink conjecture on spectral pairs\1} which encode also the information of weights, see \cite{St1}. In the simplicial case, no counter-example has been found, although it does not hold in general, see \cite{Da}. A modified version of Steenbrink's conjecture on spectral pairs has been proved in the simplicial case quite recently as a corollary of the {\it descent theorem of nearby cycle formula,} see \cite[Theorem\,\,2]{des} and (\hl{3.2.2}{3.2.2}) below. Combining this with Theorem\,\,\hl{T1}{1}, we can get the following.
\par\htt{T2}{}\msn
{\bf Theorem\,\,2.} {\it With the notation and hypotheses of Theorem\,\,$\hl{T1}{1}$, we have the equality}
\htt{8}{}
$$\aligned\SpW_f(t,u)&=\msum_{\si\ins\CFfin}\,\bl(\msum_{j=0}^{2-d_{\si}}\,t^ju^{2j}\br)q_{\si}(t)u^{d_{\si}}\\&\q{}+\bl|\CFfin^0\br|\,\bl(1\pl tu^2\br)tu\\&\q{}+{}\msum_{\si\in\CFf^0}\,(\ga_{\si}{-}3)\1t\1q_{\si}(t)\1u^2.\endaligned
\leqno(8)$$
\par\nin \ms
The right-hand side of (\hl{8}{8}) is compatible with the {\it monodromical property\1} of the weight filtration $W$, see (\hl{1.1.8}{1.1.8}) below. Note that the center of symmetry for weights is {\it always\1} $n{-}1$ independently of monodromy eigenvalues so that we have the {\it symmetry}
$$\SpW_f(t^{-1},u^{-1})\1t^nu^{2n-2}=\SpW_f(t,u).$$
\par\nin Note that the generating polynomial of the spectral pairs in the sense of \cite{SSS} is given by
$$\SpW_f(t,u^{-1})\1t^{-1}u^{2n-2}\q\h{(see also \hl{3.2}{3.2} below).}$$
\par\nin \sk
In the case $f$ has {\it non-isolated singularities,} set
$$\Sp'_f(t)\defs\Sp_f(t^{-1})\1t^n.$$
\par\nin This is called the {\it Hodge spectrum,} see \cite{GLM}, \cite{ste} and also \hl{1.1}{1.1} below.
Using the {\it Yomdin-Steenbrink formula\1} (see \hl{1.3}{1.3} below), we can deduce from Theorem\,\,\hl{T1}{1} the following.
\par\htt{T3}{}\msn
{\bf Theorem\,\,3.} {\it Assume $n\eq 3$, $f$ is Newton non-degenerate, $\Gp(f)$ is simplicial and intersects any coordinate plane of $\R^3$. For $\si\slt\Gp(f)$, let $m_{\si}$ be the number of $(d_{\si}{+}1)$-dimensional faces of $\Gp(f{+}\ell\1^r)$ $(r\,{\gg}\,0)$ which is the convex hull of $\si\cup\{r\1\ee_i\}$ for some $i\ins[1,3]$ and is not contained in any coordinate plane. Here $\ee_i$ is the \h{$i$\1th} unit vector, and $\ell$ is a sufficiently general linear function. Let $\gat_{\si}$ be the number of $2$-dimensional $($not necessarily \h{\it compact\1$)$} faces of $\Gp(f{+}\ell\1^r)$ $(r\gg 0)$ containing $\si\ins\CFf^0$. Then}
\htt{9}{}
$$\aligned\Sp'_f(t)&=\msum_{\si\in\CFfin}\,\bl(\msum_{j=0}^{2-d_{\si}}\,t^j\br)q_{\si}(t)+\bl|\CFfin^0\br|\,\bl(t{+}t^2\br)\\&\q{}+{}\msum_{\si\in\CFf^0}\,(\gat_{\si}{-}3)\1t\1q_{\si}(t)\raise12pt\h{}\\&\q{}-\msum_{\si\in\CFf}\,m_{\si}\bl(\msum_{j=0}^{1-d_{\si}}\,t^j\br)q_{\si}(t)-\msum_{\si\in\CFf^0}\,m_{\si}\1t.\endaligned
\leqno(9)$$
\par\nin \ms
Note that the multiplicities $m_{\si}$ are at most 2. Theorem\,\,\hl{T3}{3} follows from Theorem\,\,\hl{T1}{1} together with the Yomdin-Steenbrink formula for spectrum (conjectured in \cite{St2} and proved in \cite{ste}), see \hl{1.3}{1.3} and \hl{3.4}{3.4} below. Note that the latter formula does {\it not\1} hold for {\it spectral pairs,} see Remark\,\,\hl{R3.4}{3.4} below (and also \cite{KL1}, \cite{KL2}). An error in a previous version (where $m_{\si}$ did not appear) is detected during a verification of the compatibility of a recently found formula \cite[Corollary~1]{des} with Theorem\,\,\hl{T3}{3} calculating the combinatorial polynomials $r_{\!\si}(t)$. One cannot exclude the case $m_{\si}\eq2$ by assuming simply that $\Gp(f)$ intersects any coordinate plane of $\R^3$ outside the coordinate axes (for instance, $f\eq x^3y^3\pl x^2y^2z\pl xz^3\pl yz^3$).
\sk
It is unclear whether there is a simple extension of Theorems~\hl{T3}{3} to the case $n\sgt3$ (since it is not easy to see the Newton polytope). As for Theorems~\hl{T1}{1} and \hl{T2}{2}, there are rather direct generalizations as below (using the symmetry of combinatorial polynomials (\hl{6}{6}) and the estimates of the degree of $r^{\rm Def}_{\!\si}(t)$ in (\hl{7}{7})):
\par\htt{T4}{}\msn\vbox{\nin
{\bf Theorem\,\,4.} {\it Assume $n\eq 4$, and $f$ is non-degenerate with simplicial convenient Newton polytope. Let $n^{j,k}$ and $n^{j,k}_{\si}$ be the number of $\tau\ins\CFf^j$ such that $k(\tau)\eq k$ and $k(\tau)\eq k$ with $\tau\gess\si$ respectively. Then we have the equality
\htt{10}{}
$$\aligned\SpDef_f(t)&=\msum_{\si\in\CFf^0}\,\ep_{\si}\bl(t{+}t^2\br)q_{\si}(t)\\ &\q+\msum_{\tau\in\CFf^1}\,\ep_{\tau}\1 t\1q_{\tau}(t)+\ep_{\emptyset}\1t^2,\endaligned
\leqno(10)$$
\par\nin \vskip-1mm\nin
with
\vskip-7mm
\htt{11}{}
$$\aligned\ep_{\si}&:=\begin{cases}n^{1,4}_{\si}\mi 4&(k(\si)\eq 4\1),\\ n^{1,4}_{\si}\mi 1&(k(\si)\eq 3\1),\\ n^{1,4}_{\si}&(k(\si)\less2\1),\end{cases}\q\ep_{\tau}:=\begin{cases}n^{2,4}_{\tau}\mi 3&(k(\tau)\eq 4\1),\\ n^{2,4}_{\tau}\mi 1&(k(\tau)\eq 3\1),\\ n^{2,4}_{\tau}&(k(\tau)\eq 2),\end{cases}\\ \ep_{\emptyset}&:=n^{1,4}\mi 4\1n^{0,4}\mi n^{0,3}.\raise5mm\h{}\endaligned
\leqno(11)$$
\par\nin Moreover, $\ep_{\tau}\eq\ga_{\tau}{-}3$ with $\ga_{\tau}$ the number of $3$-dimensional faces of $\Gp(f)$ containing $\tau\ins\CFf^1$.}}
\par\htt{T5}{}\msn\vbox{\nin
{\bf Theorem\,\,5.} {\it With the notation and assumptions of Theorem\,\,$\hl{T4}{4}$, we have the equality
\htt{12}{}
$$\aligned\SpDefW_f(t,u)&=\msum_{\si\in\CFf^0}\,\ep_{\si}\bl(tu^2{+}t^2u^4\br)q_{\si}(t)\raise15pt\h{}\\ &\q+\msum_{\tau\in\CFf^1}\,\ep_{\tau}\1 t\1q_{\tau}(t)u^3+\ep_{\emptyset}\1t^2u^3.\endaligned
\leqno(12)$$
\par\nin Here $\,\SpDefW_f(t,u)$ is the weighted defect as in {\rm\hl{3.2}{3.2}} below.}}
\ms
These follow from (\hl{5}{5}) (using \cite{exp}) and \cite[Theorem 2]{des} together with (\hl{6}{6}) (that is, Theorem\,\,\hl{TA}{A} in Appendix) and (\hl{7}{7}), see \hl{3.3}{3.3} below.
The formulas become more and more complicated as the dimension $n$ increases.
Note that a Thom-Sebastiani type formula does not hold for $\Ga$-spectrum nor for defective spectral numbers (for instance, if $f\eq g(x,y)\pl z^c$ with $g$ non-degenerate).
In general (except for $n\eq 3$), it is unclear whether the coefficients of $\SpDef_f(t)$ are {\it non-negative,} see Remark\,\,\hl{RA}{A} in Appendix below.
\sk
The formula in Theorem\,\,\hl{T1}{1} was first found as a consequence of many calculations using a small computer program for $n\eq 3$.
It is quite important to find simple formulas such as Theorems~\hl{T1}{1}, \hl{T2}{2} and \hl{T3}{3} in the 3 variable case, since we can calculate explicit examples effectively only in this case. When $n\sgt3$, we may have to write a computer program determining the Newton polytopes at least. (There is an interesting quiz: what is the minimal triple of positive integers $(a,b,c)$ in lexicographic order satisfying the following condition? There is a polynomial $f$ of 4 variables which is invariant by the action of the symmetric group ${\mathfrak S}_4$ and whose Newton polytope $\Gp(f)$ contains $(a,a,a,a)$, $(b,b,b,0)$, $(c,0,0,0)$ as vertices and also the convex hull of $\{(a,a,a,a),(b,b,b,0),(b,b,0,b),(b,0,b,b)\}$ as a face. The answer may be much larger than is expected, even in the corresponding case for $n\eq 3$, see \hl{2.4}{2.4} below.)
\sk
In Section~1 we review some basics of spectrum, Newton filtration, and the Yomdin-Steenbrink formula. In Section~2 we present a precise proof of Steenbrink's formula for the spectrum in the convenient non-degenerate case. In Section~3 we prove the main theorems after introducing $\Ga$-spectrum and explaining the Steenbrink conjecture on spectral pairs.
\bs\bs
\centerline{Contents}
\ms
\par\hl{S1}{1. Preliminaries}\hfill 6
\par\q\hl{1.1}{1.1.~Spectrum}\hfill 6
\par\q\hl{1.2}{1.2.~Non-degenerate case}\hfill 7
\par\q\hl{1.3}{1.3.~Yomdin-Steenbrink formula}\hfill 8
\par\q\hl{1.4}{1.4.~Non-isolated non-degenerate case}\hfill 9
\par\q\hl{1.5}{1.5.~Monodromy zeta functions}\hfill 13
\par\hl{S2}{2. Proof of Steenbrink's formula}\hfill 13
\par\q\hl{2.1}{2.1.~Steenbrink formula for spectrum}\hfill 13
\par\q\hl{2.2}{2.2.~Determination of Newton polyhedra for $n\eq 3$}\hfill 18
\par\q\hl{2.3}{2.3.~Topological classification of Newton boundaries for $|\CFf^0|\eq6$, $n\eq3$}\hfill 19
\par\q\hl{2.4}{2.4.~Examples for $n\eq 4$}\hfill 20
\par\hl{S3}{3. Simplicial Newton polytope case}\hfill 22
\par\q\hl{3.1}{3.1.~$\Ga$-spectrum}\hfill 22
\par\q\hl{3.2}{3.2.~Spectral pairs}\hfill 24
\par\q\hl{3.3}{3.3.~Proofs of Theorems}\,\,\hl{T1}{1}--\hl{T2}{2} and \hl{T4}{4}--\hl{T5}{5}\hfill 26
\par\q\hl{3.4}{3.4.~Proof of Theorem}\,\,\hl{T3}{3}\hfill 26
\par\hl{A}{Appendix. Symmetry of combinatorial polynomials}\hfill 29
\par\bs\bs\htt{S1}{}
\vbox{\centerline{\bf 1. Preliminaries}
\bsn
In this section we review some basics of spectrum, Newton filtration, and the Yomdin-Steenbrink formula.}
\par\htt{1.1}{}\msn
{\bf 1.1.~Spectrum.} Let $f\,{:}\,(\C^n,0)\tos(\C,0)$ be a holomorphic function. We have a canonical mixed Hodge structure on the (reduced) vanishing cohomology $\Ht^j(\Ff,\Q)$, where $\Ff$ denotes the Milnor fiber of $f$. The $\la$-eigenspace of the semisimple part $T_s$ of the {\it Jordan decomposition\1} of the monodromy $T$ is denoted by $\Ht^j(\Ff,\C)_{\la}$. (Note that $T$ is the {\it inverse\1} of the Milnor monodromy, see for instance \cite{DS}.)
\sk
The {\it Steenbrink spectrum\1} $\Sp_f(t)\eq\sum_{\al\in\Q_{>0}}n_{\al}t^{\1\al}$ (see \cite{St1}, \cite{St2}) is a fractional power polynomial defined by
\htt{1.1.1}{}
$$\aligned&n_{\al}=\msum_{j=0}^{n-1}\,(-1)^j\dim_{\C}\Gr_F^p\Ht^{n-1-j}(\Ff,\C)_{\la}\\&\q\q\h{for}\q\,\,p\eq[n{-}\al],\q\la\eq\exp(-2\pi i\al).\endaligned
\leqno(1.1.1)$$
\par\nin \sk
The {\it Hodge spectrum\1} $\Sp'_f(t)$ (see \cite{ste}, \cite{GLM}) is the fractional power polynomial
\htt{1.1.2}{}
$$\Sp'_f(t)\defs\Sp_f(t^{-1})\1t^n.
\leqno(1.1.2)$$
\par\nin This can be defined also by $\Sp'_f(t)\eq\sum_{\al\in\Q_{>0}}n'_{\al}t^{\1\al}$ with
\htt{1.1.3}{}
$$\aligned&n'_{\al}=\msum_{j=0}^{n-1}\,(-1)^j\dim_{\C}\Gr_F^p\Ht^{n-1-j}(\Ff,\C)_{\la}\\&\q\q\h{for}\q\,\,p\eq[\al],\q\la\eq\exp(2\pi i\al).\endaligned
\leqno(1.1.3)$$
\par\nin \sk
In the isolated singularity case, the $n_{\al}$ are non-negative, and we have
\htt{1.1.4}{}
$$\Sp_f(t)\eq\msum_{i=1}^{\mu_f}\,t\1^{\al_i},
\leqno(1.1.4)$$
\par\nin where $\mu_f$ is the Milnor number of $f$. The $\al_i$ are positive rational numbers which are assumed to be weakly increasing, and are called the {\it spectral numbers\1} of $f$. We have the {\it symmetry}
\htt{1.1.5}{}
$$\al_i\pl\al_j\eq n\,\,\,\,(\1 i{+}j\eq\mu_f{+}1\1),\q\h{that is,}\q\Sp_f(t)\eq\Sp'_f(t).
\leqno(1.1.5)$$
\par\nin This follows from {\it Hodge decomposition\1} of the graded quotients of the weight filtration $W$ together with the {\it monodromical property\1} of $W$ as in (\hl{1.1.8}{1.1.8}) below (using also the assertion that the monodromy is defined on $H^{\ssb}(\Ff,\R)$), see for instance \cite{St1}.
\sk
In the isolated hypersurface singularity case, we have
\htt{1.1.6}{}
$$n_{\al}\eq\dim\Gr_V^{\al}\bl(\C\{x\}/(\dd f)\br)\q(\al\ins\Q).
\leqno(1.1.6)$$
\par\nin Here $(\dd f)\sst\C\{x\}\eq\OO_{\C^n,0}$ denotes the Jacobian ideal with $(x)\eq(x_1,\dots,x_n)$ the canonical coordinate system of $\C^n$, and $V$ on $\C\{x\}/(\dd f)$ is induced from the $V$-filtration on the {\it Brieskorn lattice\1} $H''_f\defs\Om_{\C^n,0}^n/\ddd f{\wedge}\1\ddd\Om_{\C^n,0}^{n-2}$. Moreover there are canonical isomorphisms
\htt{1.1.7}{}
$$\aligned&\Gr_V^{\al}\bl(\C\{x\}/(\dd f)\br)=\Gr_F^pH^{n-1}(\Ff,\C)_{\la}\\ &\q\q\bl(p\eq[n{-}\al],\,\,\,\,\la\eq\exp(-2\pi i\al)\br).\endaligned
\leqno(1.1.7)$$
\par\nin Note that the graded action of $\Gr_Vf$ on the left-hand side (with index $\al$ shifted by 1) corresponds to that of $-N/2\pi i$ on the right-hand side (with index $p$ shifted by $-1$) via these isomorphisms, where $N\defs\log T_u$ with $T_u$ the unipotent part of the {\it Jordan decomposition\1} $T\eq T_sT_u$ of the monodromy, see for instance \cite{SS}.
\par\htt{R1.1a}{}\msn
{\bf Remark\,\,1.1a.} The spectrum in \cite{St2}, \cite{Sing} is shifted by $-1$ compared with the above definition which agrees with \cite{ste}, \cite{St1}, etc.
\par\htt{R1.1b}{}\msn
{\bf Remark\,\,1.1b.} The difference between $\Sp_f(t)$ and $\Sp'_f(t)$ is closely related to the one between the cohomological pull-back functors $H^ki_0^*$ and $H^{-k}i_0^!$ of mixed Hodge modules, which are dual to each other up to a Tate twist depending on monodromy eigenvalues, where $i_0:\{0\}\into\C^n$ is the inclusion, see for instance \cite[(2.6.2)]{mhm}.
\par\htt{R1.1c}{}\msn
{\bf Remark\,\,1.1c.} In the isolated singularity case, we have to use the decomposition by the
{\it unipotent\1} and {\it non-unipotent\1} monodromy part of the vanishing cohomology:
$$H^{n-1}(\Ff,\C)\eq H^{n-1}(\Ff,\C)_1\oplus H^{n-1}(\Ff,\C)_{\ne 1},$$
\par\nin where the last term is defined by
$$H^{n-1}(\Ff,\C)_{\ne 1}\defs\mopl_{\la\ne 1}\,H^{n-1}(\Ff,\C)_{\la}.$$
\par\nin The weight filtration $W$ on $H^{n-1}(\Ff,\C)_{\ne 1}$ and $H^{n-1}(\Ff,\C)_1$ is given by the {\it monodromy filtration\1} shifted respectively by $n{-}1$ and $n$, that is,
\htt{1.1.8}{}
$$\aligned N^j:\Gr_{n-1+j}^WH^{n-1}(\Ff,\C)_{\ne 1}&\simto\Gr_{n-1-j}^WH^{n-1}(\Ff,\C)_{\ne 1}\q(j\ins\N),\\ N^j:\Gr_{n+j}^WH^{n-1}(\Ff,\C)_1&\simto\Gr_{n-j}^WH^{n-1}(\Ff,\C)_1\q(j\ins\N),\endaligned
\leqno(1.1.8)$$
\par\nin where $N\eq\log T_u$ as in a remark after (\hl{1.1.7}{1.1.7}). This implies that the vanishing cohomology $H^{n-1}(\Ff,\C)_{\ne 1}$, $H^{n-1}(\Ff,\C)_1$ have pure weight $n{-}1$ and $n$ respectively if $N\eq0$ (for instance, if $f$ is weighted homogeneous, where the Milnor fiber $\Ff$ is algebraic). The proof of (\hl{1.1.8}{1.1.8}) is not quite trivial, since some argument as in the proof of \cite[4.2.2]{mhp} is required. (Indeed, the weight filtration $W$ is {\it not\1} defined as the shifted monodromy filtration in \cite{St1}, and the passage from the $E_1$-term to the $E_2$-term is highly non-trivial, see also \cite{GN}.) Note that $N$ is a morphism of mixed Hodge structures of type $(-1,-1)$ decreasing the indices of $F,W$ by 1 and 2 respectively, and is bistrictly compatible with a pair of two filtrations $(F,W)$ so that the kernel and cokernel of $N^i$ commute with the passage to the graded quotients $\Gr^{\ssb}_F\Gr_{\ssb}^W$.
\par\htt{1.2}{}\msn
{\bf 1.2.~Non-degenerate case.} We denote by $\Gp(f)$ the {\it Newton polytope\1} of $f\ins\C\{x\}$. This is by definition the convex hull of the union of $\nu\pl\R_{\ges 0}^n$ for $\nu\ins\Spx\1f$, where
\htt{1.2.1}{}
$$\Spx\1f:=\{\nu\ins\N^n\mid c_{\nu}\ne 0\1\}\q\h{for}\q f\eq\msum_{\nu}\,c_{\nu}x^{\nu}\ins\C\{x\},
\leqno(1.2.1)$$
\par\nin with $x_1,\dots,x_n$ the coordinates of $\C^n$.
\sk
We say that $f$ is ({\it Newton\1}) {\it non-degenerate,} or more precisely, $f$ has {\it non-degenerate Newton boundary} (or {\it principal} \h{\it part\1}) if for any {\it compact\1} face $\si\sst\Gp(f)$, we have
\htt{1.2.2}{}
$$\mcap_{i=1}^n\,\bl\{\1 x_i\dd_{x_i}f_{\si}\eq0\br\}\sst\{\1x_1\cdots x_n\eq0\1\},
\leqno(1.2.2)$$
\par\nin where $f_{\si}\defs\mopl_{\nu\in\si}\,c_{\nu}x^{\nu}$ with $c_{\nu}$ as in (\hl{1.2.1}{1.2.1}), see \cite{Ko}, \cite{Va1}, etc.
\sk
Assume $f$ is non-degenerate, and is {\it convenient,} that is, its Newton polytope $\Gp(f)$ is. (The latter condition says that $\Gp(f)$ intersects every coordinate axis of $\R^n$.) Note that these conditions imply that $\mu_f\eq\dim\C\{x\}/(\dd f)\slt\infty$, that is, $f$ has an isolated singularity at 0, see \cite{Ko}. For $g\ins\C\{x\}$, set
\htt{1.2.3}{}
$$v_f(g):=\max\bl\{r\ins\R\mid\ob\pl\Spx\1g\,\subset\,r\,\Gp(f)\br\},
\leqno(1.2.3)$$
\par\nin with $\ob\defs(1,\dots,1)$. We have $\al_1\eq v_f(1)$, and the $V$-filtration on $\C\{x\}/(\dd f)$ in (\hl{1.1.6}{1.1.6}) is induced from the Newton filtration $V_N^{\ssb}$ on $\C\{x\}$ defined by
\htt{1.2.4}{}
$$V_N^{\be}\C\{x\}:=\bl\{g\ins\C\{x\}\mid v_f(g)\gess\be\br\},
\leqno(1.2.4)$$
\par\nin see \cite{exp}.
\sk
For {\it logarithmic\1} complexes having logarithmic poles along coordinate hyperplanes, we need the {\it unshifted\1} Newton filtration $V'_N$ on $\C\{x\}$ which is defined by omitting ``$\,\ob+\,$" in (\hl{1.2.3}{1.2.3}--\hl{1.2.4}{4}).
\par\htt{R1.2a}{}\msn
{\bf Remark\,\,1.2a.} The convenience condition can be {\it avoided\1} in the main theorem of \cite{exp} as is shown in \cite{JKSY}, assuming the singularity is isolated. Note also that we can assume $f$ is convenient in the isolated hypersurface singularity case by adding monomials $x_i^{a_i}$ to $f$ for $a_i\gg 0$ if we use finite determinacy of holomorphic functions with isolated singularities.
(In \cite{VK}, {\it another\1} proof of the main theorem of \cite{exp} is given, where the {\it symmetry\1} of the Newton spectrum seems to be shown in a rather technical way; for instance, using the relation between the Laurent expansions of a certain rational function of $t$ at 0 and at $\infty$. Note that it is rather easy to show the inclusion $V_N^{\ssb}\sst V^{\ssb}$ on $\C\{x\}/(\dd f)$.)
\par\htt{1.2b}{}\msn
{\bf Remark\,\,1.2b.} In our case condition (\hl{1.2.2}{1.2.2}) for non-degenerate functions is equivalent to the {\it smoothness\1} of the hypersurface in $(\C^*)^n$ defined by $f_{\si}$, since $f_{\si}$ is a linear combination of $x_i\dd_{x_i}f_{\si}$ ($i\ins[1,n]$) using the {\it positivity\1} of the coefficients of linear functions defining compact faces of $\Gp(f)$, see \cite{JKSY}. In the Laurent polynomial case this does not hold, if we define the Newton boundary {\it naively\1} without taking care of the origin as in \cite[Definition~1.5]{Ko}; for instance, if $f\eq x_1(x_1\pl1)^2+\msum_{i=2}^n\,x_i^{a_i}\pl c\,$ with $a_i\gess 2$ and $c\ins\C^*$ sufficiently general, then the smoothness holds, but (\hl{1.2.2}{1.2.2}) fails, when $\si$ is the face contained in the $x_1$-axis. This shows that smoothness may be insufficient in certain cases.
\par\htt{1.3}{}\msn
{\bf 1.3.~Yomdin-Steenbrink formula.} Using the {\it Hodge spectrum\1} $\Sp'_f(t)$ instead of $\Sp_f(t)$ (see \hl{1.1}{1.1}), the {\it Yomdin-Steenbrink formula\1} can be described as follows:
\sk
Assume for simplicity the singular locus of $(f^{-1}(0),0)$ is a union of {\it smooth\1} curves $(C_i,0)$ for $i\ins I$. On $C_i^*\defs C_i\stm\{0\}$, the vanishing cohomology sheaf of $f$ underlies a variation of mixed Hodge structure. Its stalk coincides with the vanishing cohomology of the restriction $f_{(i)}$ of $f$ to a hyperplane transversal to $C_i^*$. Let $\al_{i,j}$ be its spectral numbers ($j\ins[1,\mu_{f_{(i)}}])$. Since the monodromy of local system on $C_i^*$ commutes with the Milnor monodromy, there are rational numbers $\be_{i,j}\ins[0,1)$ such that $e^{2\pi\sqrt{-1}\1\be_{i,j}}$ is the eigenvalue of the local system monodromy for an eigenvector associated with $\al_{i,j}$. Here the well-definedness of $\be_{i,j}$ does not seem very clear unless the Hodge filtration on the graded pieces of the weight filtration on the vanishing cohomology sheaf at general points of the singular locus of $f$ is locally constant, that is, the filtration is stable by the action of vector fields. It is more precise to say that the $\be_{i,j}$ are determined in the proof of the formula. (Note that $[0,1)$ and $e^{2\pi\sqrt{-1}\1\be_{i,j}}$ are respectively replaced by $(0,1]$ and $e^{-2\pi\sqrt{-1}\1\be_{i,j}}$ if $\Sp_f(t)$ is used.)
\sk
Let $\ell$ be a sufficiently general linear function (fixing some local coordinates) such that $\ell^{-1}(0)$ is transversal to any $C_i$, and is not contained in the set of limits of tangent spaces of smooth points of $f^{-1}(0)$ (which has dimension at most $n-2$). Then the restriction of $f$ to $\ell^{-1}(0)$ has an isolated singularity at $0$, and so does $f\pl\ell\1^r$ for $r\gg 0$. We have moreover the following equalities (see \cite{ste}):
\htt{1.3.1}{}
$$\Sp_{f+\ell\1^r}(t)-\Sp'_f(t)=\msum_{i,j}\,\msum_{k=0}^{r-1}\,t\1^{\al_{i,j}\1+\1\be_{i,j}/r\1+\1 k/r}\q(r\gg 0).
\leqno(1.3.1)$$
\par\nin This is the Yomdin-Steenbrink formula, which is conjectured in \cite{St2}, and proved in \cite{ste}.
Note that we have $\Sp_{f+\ell\1^r}(t)\eq\Sp'_{f+\ell\1^r}(t)$.
We can then get the Hodge spectrum $\Sp'_f(t)$, if we know the $\al_{i,j}$ and if the $\be_{i,j}$ can be determined by calculating the differences
\htt{1.3.2}{}
$$\aligned&\Sp_{f+\ell\1^r}(t)-\Sp_{f+\ell\1^{r'}}(t)=\Sp_{f+\ell\1^r}(t)''-\Sp_{f+\ell\1^{r'}}(t)'',\\&\h{with}\q\q\Sp_{f+\ell\1^r}(t)''\defs\msum_{i,j}\,\msum_{k=0}^{r-1}\,t\1^{\al_{i,j}\1+\1\be_{i,j}/r\1+\1 k/r},\endaligned
\leqno(1.3.2)$$
\par\nin for sufficiently large two distinct prime numbers $r,r'$.
\sk
It might be expected that one could get a formula for the Hodge spectrum $\Sp'_f(t)$ using (\hl{1.3.1}{1.3.1}--\hl{1.3.2}{2}). Unfortunately, this is not necessarily easy without combining it with some other method. One possibility is explained in \hl{1.4}{1.4} below.
\par\htt{R1.3}{}\msn
{\bf Remark\,\,1.3.} In the case $n\eq 3$ we may then assume that
\htt{1.3.3}{}
$$\al_{i,j}\pl\al_{i,j'}\eq2,\,\,\,\,\be_{i,j}\pl\be_{i,j'}\ins\Z\,\,\,\,\,\h{for}\,\,\, j{+}j'\eq\mu_{(i)}{+}1,
\leqno(1.3.3)$$
\par\nin with $\mu_{(i)}$ the Milnor number of $f_{(i)}$, renumbering the $\al_{i,j},\be_{i,j}$ so that the eigenvector for $\al_{i,j},\be_{i,j}$ and the one for $\al_{i,j'},\be_{i,j'}$ are complex conjugate to each other (since the complex conjugate of $e^{2\pi i\al}$ is $e^{-2\pi i\al}$ for $\al\ins\R$).
\sk
Indeed, the monodromies are defined over $\R$, that is, on the cohomology with $\R$-coefficients, and they are compatible with the associated complex conjugation. Here we use also the Hodge decomposition of the non-unipotent monodromy part of the vanishing cohomology of $f_{(i)}$, which has pure weight 1. (This can be generalized to the case $n\sgt3$ by the $N$-primitive decomposition of the graded quotients of the weight filtration.)
\sk
Using the above construction, Question~\hl{Q1.4}{1.4} below for $n\eq 3$ can be reduced essentially to the case $\al_{i,j}\less 1$. Indeed, the right-hand side of (\hl{1.3.1}{1.3.1}) is determined by the $\al_{i,j},\be_{i,j}$ with $j\less(\mu_{(i)}\pl1)/2$ by (\hl{1.3.3}{1.3.3}).
\par\htt{1.4}{}\msn
{\bf 1.4.~Non-isolated non-degenerate case.} It is not necessarily easy to {\it find\1} a simple formula for the spectrum $\Sp(t)$ or $\Sp'(t)$ in the {\it non-degenerate non-isolated singularity\1} case (for instance, the calculation of combinatorial polynomials is rather complicated in \cite[Corollary~1]{des}). This is quite different from the case of monodromy zeta functions where we have a simple formula using $\Gp(f)$, see \cite{Va1}. Note that the combinatorial polynomial is simplified very much for these. In this section we show one possibility using the Yomdin-Steenbrink formula (\hl{1.3.1}{1.3.1}).
\sk
Set
$$\aligned I_f&:=\bl\{i\ins[1,n]\mid\Spx f\cap\N\1\ee_i\eq\emptyset\br\},\\
f_{(i)}&:=f|_{x_i=\ga_i}\,\,\,\h{for $\ga_i\ins\De_{\ep}^*$ sufficiently general}\,\,(i\ins I_f).\endaligned$$
\par\nin Here $\ee_i$ is the \h{$i$\1th} unit vector, and $\De_{\ep}^*$ is a punctured disk of radius $\ep\ll 1$.
\sk
Assume $I_f\ne\emptyset$ and the following:
\htt{SC}{}
$$\aligned&\h{$f\,$ has simplicial non-degenerate Newton boundary, and the $f_{(i)}$}\\&\h{are convenient (with $\ga_i\ins\C$ sufficiently general) for any $i\ins I_f$.}\endaligned
\leqno{\rm(SC)}$$
\par\nin \sk
By \cite[Proposition\,\,A.2]{JKSY}, we see that $f\pl\ell\1^r$ is {\it non-degenerate\1} for $r$ sufficiently large. (Indeed, the convenience of $f_{(i)}$ implies that $\Spx f\cap\bl(\N\ee_i\pl\N\ee_j)\ne\emptyset$ for any $i,j\ins[1,n]$ with $i\ne j$. Hence $\{r\1\ee_i,r\1\ee_j\}$ cannot be contained in the same compact face of $\Gp(f{+}\ell\1^r)$ when $r\gg 0$.) So the spectrum $\Sp_{f+\ell\1^r}(t)$ can be computed using the Newton filtration $V_N^{\ssb}$ on $\C\{x\}/(\dd f)$ ($r\gg 0$), see \cite{exp}.
\sk
The non-degeneracy of $f_{(i)}$ in (\hl{SC}{SC}) seems to follow from that of $f$ if $n\eq 3$ (using a theory on equisingularities of plane curves). This is, however, unnecessary for us using Lemma\,\,\hl{L1.4}{1.4} below.
\sk
For $\,\nu^{(i)}\eq\bl(\nu_k^{(i)}\br)_{k\ne i}\ins\N^{(i),n-1}\,\,(i\ins I_f$) with $\N^{(i),n-1}\defs\{\nu\ins\N^n\mid\nu_i\eq0\}$, set
$$x^{\nu^{(i)}}:=\mprod_{k\ne i}\,x_k^{\nu_k^{(i)}}\ins\C\{x_{(i)}\},$$
\par\nin where $x_{(i)}\eq(x_k)_{k\ne i}$ is the coordinate system of the target of the projection $\pi_i:\C^n\tos\C^{n-1}$ forgetting the \h{$i$\1th} coordinate, see also \cite[(A.1.4)]{JKSY}.
\sk
In order to make use of the Yomdin-Steenbrink formula (\hl{1.3.1}{1.3.1}), we have to lift some monomials in $\C\{x_{(i)}\}$ to certain monomials in $\C\{x\}$. We have the following.
\par\htt{Q1.4}{}\msn
{\bf Question\,\,1.4.} Are there $\,\nu^{(i),j}\eq\bl(\nu_k^{(i),j}\br)_{k\ne i}\ins\N^{(i),n-1}\,\,\,\,(i\ins I_f,\,j\ins[1,\mu_{f_{(i)}}])\,$ satisfying the two conditions (a) and (b) below under the assumption (\hl{SC}{SC})?
\par\htt{a}{}\msn
(a) In the notation of (\hl{1.2.3}{1.2.3}--\hl{1.2.4}{4}), we have
\htt{1.4.1}{}
$$v_f\bl(x^{\nu^{(i),j}}\br)\eq\al_{i,j},\,\,\,\,\h{in particular,}\,\,\,\,x^{\nu^{(i),j}}\ins V_N^{\al_{i,j}}\C\{x_{(i)}\},
\leqno(1.4.1)$$
\par\nin and moreover the $\,\bl[x^{\nu^{(i),j}}\br]\ins\Gr_{V_N}^{\al}\C\{x_{(i)}\}$ with $\al_{i,j}\eq\al$ are linearly independent over $\C$ for any $\al\ins\Q$.
(The last condition is actually equivalent to that the $x^{\nu^{(i),j}}$ with $\al_{i,j}\eq\al$ are mutually distinct, since the filtration $V_N$ has a canonical splitting by monomials. It may be expected that this condition holds even modulo the Jacobian ideal $(\dd f_{(i)})$, but this stronger assumption is not needed for our purpose. Note also that this can be satisfied if we have the equality $\Sp_{f_{(i)}}(t)\eq\SpG_{f_{(i)}}(t)$ for any $i\ins I_f$, in particular, if $n\less 3$.)
\par\htt{b}{}\msn
(b) Let $\tau^{(i),j}$ be an $(n{-}2)$-dimensional compact face of $\Gp(f_{(i)})$ such that $\nu^{(i),j}\ins\N^{(i),n-1}$ is contained in its cone. Let $\si^{(i),j}$ be an $(n{-}2)$-dimensional compact face of $\Gp(f)$ such that its image $\pi_i(\si^{(i),j})$ is contained in $\tau^{(i),j}$ and the cone of the image contains $\nu^{(i),j}$. Let $\nu_i^{(i),j}$ be a positive rational number satisfying
\htt{1.4.2}{}
$$\nut^{(i),j}:=\bl(\nu_k^{(i),j}\br)_{k\in[1,n]}\ins C(\si^{(i),j}),
\leqno(1.4.2)$$
\par\nin where $C(\si^{(i),j})$ is the cone of $\si^{(i),j}$ in the vector space $\R^n$. Then we have the equalities
\htt{1.4.3}{}
$$\be_{i,j}=\bl\{-\nu_i^{(i),j}\br\}\q\bl(i\ins I_f,\,j\ins[1,\mu_{f_{(i)}}]\br),
\leqno(1.4.3)$$
\par\nin where $\{\al\}\defs\al-[\al]\ins[0,1)$ for $\al\ins\R$.
\ms
We may consider the restriction of Question~\hl{Q1.4}{1.4} to the case $\al_{i,j}\less 1$ (especially in the case $n\eq 3$ using Remark\,\,\hl{R1.3}{1.3}).
\sk
We see that Question~\hl{Q1.4}{1.4} has a positive answer in many cases as follows.
\par\htt{E1.4a}{}\msn
{\bf Example~1.4a.} Let $f\eq x^ay^b$ ($a,b\gess 2$). Set $e\defs{\rm GCD}(a,b)$. For $r\sgt a\pl b$, it is well-known (see \cite{St1}, \cite{Ar}) that
$$\aligned\Sp'_f(t)&=-\1\msum_{j=1}^{e-1}\,t\1^{j/e}+\msum_{j=e}^{2e-1}\,t\1^{j/e},\\ \Sp_{f+\ell\1^r}(t)&=\Sp'_f(t)+\msum_{j=1}^{a-1}\,\msum_{k=0}^{r-1}\,t\1^{j/a+(\{-jb/a\}\pl k)/r}\\&\q\q\q\q+\msum_{j=1}^{b-1}\,\msum_{k=0}^{r-1}\,t\1^{j/b+(\{-ja/b\}\pl k)/r}.\endaligned$$
\par\nin We have $\al_{1,j}\eq\tfrac{j}{b}$ ($j\ins[1,b{-}1]$), $\al_{2,j}\eq\tfrac{j}{a}$ ($j\ins[1,a{-}1]$), and $\mu_{f+\ell\1^r}\eq(a{+}b{-}2)r{+}1$. Hence
$$\be_{1,j}\eq\{-ja/b\},\q\be_{2,j}\eq\{-jb/a\}.$$
\par\nin This formula for the $\be_{i,j}$ is also well-known, see for instance \cite[Theorem\,\,3.3]{mhm}. This essentially gives a general formula for the $\be_{i,j}$ in the general case with $n\eq2$.
\par\htt{E1.4b}{}\msn
{\bf Example~1.4b.} Let $f\eq x^4y^2\pl x^2y^3$. We have $\mu_{f+\ell\1^r}\eq2r\pl9$ and
$$\Sp_{f+\ell\1^r}(t)=\msum_{j\in J}\,t\1^{j/8}+2\1t+\msum_{k=1}^{r-1}\,t\1^{1/2\1+\1 k/r}+\msum_{k=0}^{r-1}\,t\1^{1/2\1+\1 1/2r\1+\1 k/r},$$
\par\nin where $J\defs\{3,4,6,7,9,10,12,13\}$ and $r\gess 9$. We have the following picture of monomials corresponding to the spectral numbers of $f\pl\ell\1^9$ as in \cite{Ar}:
\ms
$$\setlength{\unitlength}{4mm}
\begin{picture}(13,12)
\multiput(0,0)(1,0){14}{\line(0,1){12}}
\multiput(0,0)(0,1){13}{\line(1,0){13}}
\multiput(1,1)(1,0){10}{\circle*{.35}}
\multiput(1,4)(0,1){7}{\circle*{.35}}
\multiput(1,2)(1,0){4}{\circle*{.35}}
\multiput(1,3)(1,0){4}{\circle*{.35}}
\put(5,4){\circle*{.35}}
\put(6,3){\circle*{.35}}
\put(2,1){\circle{.5}}
\put(0,9){\line(1,-3){2}}
\put(2,3){\line(2,-1){2}}
\put(4,2){\line(5,-2){5}}
\thicklines
\put(0,0){\line(2,1){8}}
\put(2,3){\line(2,1){4}}
\put(9,0){\line(2,1){4}}
\put(0,0){\line(1,0){9}}
\put(4,2){\line(1,0){9}}
\put(0,0){\line(2,3){4}}
\put(4,2){\line(2,3){2}}
\put(0,9){\line(2,3){2}}
\put(0,0){\line(0,1){9}}
\put(2,3){\line(0,1){9}}
\end{picture}$$
\par\nin \msn
Here the black vertices give the monomials corresponding to spectral numbers (see (\hl{3.1.1}{3.1.1}) below), and the vertex at $(2,1)$ corresponds to the difference between $J$ and $J'$ explained just below.
We see that $\mu_{f_{(k)}}\eq1$, $\al_{k,1}\eq\tfrac{1}{2}$ ($k\eq1,2$), and hence
$$\be_{1,1}\eq0,\q\be_{2,1}\eq\tfrac{1}{2},\q\Sp'_f(t)\eq\msum_{j\in J'}\,t\1^{j/8}\pl2\1t.$$
\par\nin Here $J\eq J\stm\{4\}$, since the second summation on the right-hand side of the above equality starts with $k\eq1$. This is compatible with \cite{Va1}.
\par\htt{1.4c}{}\msn
{\bf Example~1.4c.} Let $f\eq x^2y^2\pl z^4$. This is a simple example such that we have $m_{\si}\eq2$ for some $\si\ins\CFfin^1$. Here
$$\{\al_{i,j}\}_{j\in[1,3]}=\bl\{\tfrac{3}{4},1,\tfrac{5}{4}\br\}\q(i\eq 1,2).$$
\par\nin Using \cite{exp}, we can verify that
$$\aligned\Sp_{f+\ell\1^r}(t)&=-t^{3/4}\mi t\pl t^{3/2}\pl 2\1t^{7/4}\pl t^2\pl t^{9/4}\\&{}\,+2\,\msum_{k=0}^{r-1}\,\bl(t\1^{3/4+k/r}\pl t\1^{1+k/r}\pl t\1^{5/4+k/r}\br),\endaligned$$
\par\nin with $\mu_{f+\ell\1^r}\eq6\1r\pl 3$ ($r\gg 0$). We thus get (in a compatible way with \cite{Va1})
$$\Sp'_f(t)=-t^{3/4}\mi t\pl t^{3/2}\pl 2\1t^{7/4}\pl t^2\pl t^{9/4},$$
\par\nin with $\be_{i,j}\eq0$ ($\forall\,i,j$). This agrees with the Thom-Sebastiani type theorem for spectrum (see for instance \cite{MSS}), since
$$\aligned&\bl(-t^{1/2}{+}t{+}t^{3/2}\br)\bl(t^{1/4}{+}t^{1/2}{+}t^{3/4}\br)\\ &=-t^{3/4}\mi t\pl t^{3/2}\pl 2\1t^{7/4}\pl t^2\pl t^{9/4}.\endaligned$$
\par\nin \ms
We have more generally the following.
\par\htt{P1.4}{}\msn
{\bf Proposition\,\,1.4.} {\it Question~{\rm\hl{Q1.4}{1.4}} is solvable if $n\less 3$.}
\msn
{\it Proof.} Since the argument in the case $n\eq2$ is essentially the same as in Examples~\hl{E1.4a}{1.4a}--\hl{E1.4b}{b}, we may assume $n\eq 3$.
We have $\Sp_{f_{(i)}}(t)\eq\SpG_{f_{(i)}}(t)$ ($i\ins I_f$), since $n{-}1\eq 2$. So there are $\nu^{(i),j}\ins\N^{(i),n-1}$ for $i\ins I_f$, $j\ins[1,\mu_{f_{(i)}}]$ such that (\hl{1.4.1}{1.4.1}) holds by the definition of $\Ga$-spectral numbers. This is just Arnold's picture (see also Example~\hl{E1.4b}{1.4b}), since the $f_{(i)}$ have two variables. Here it is not questioned whether the monomials $x^{\nu^{(i),j}}$ ($j\ins[1,\mu_{f_{(i)}}]$) give a basis of $\C\{x_{(i)}\}/(\dd f_{(i)})$.
\sk
For condition (\hl{b}{b}), let $\ell$ be the linear function (without a constant term) on $\R^n$ such that $\ell^{-1}(1)$ contains $\si^{(i),j}\cup\{r\1\ee_i\}$ ($r\gg 0$). Let $\ell'$ be the linear function on $\R^{n-1}$ with $\ell'{}^{-1}(1)\supset\tau^{(i),j}$. Then we have the equality
\htt{1.4.4}{}
$$\ell\bl(\nut^{(i),j}\br)\eq\ell'\bl(\nu^{(i),j}\br)\,({=}\,\al_{i,j}),
\leqno(1.4.4)$$
\par\nin since $\pi_i(\si^{(i),j})\sst\tau^{(i),j}$ and $\pi_i(\nut^{(i),j})\eq\nu^{(i),j}$. (Indeed, (\hl{1.4.4}{1.4.4}) follows by restricting the linear functions $\ell$, $\ell'$ respectively to the $1$-dimensional vector subspaces spanned by $\nut^{(i),j}$ and $\nu^{(i),j}$, which are identified to each other by the projection $\pi_i$.)
\sk
The assertion (\hl{1.4.3}{1.4.3}) then follows from the equality (\hl{1.4.4}{1.4.4}) by comparing (\hl{1.3.2}{1.3.2}) with (\hl{4}{4}) and renumbering the $\nu^{(i),j}$ using a permutation of $[1,\mu_{f_{(i)}}]$. Here we consider monomials $\,x^{\nut^{(i),j,k}}\,$ for $\,\nut^{(i),j,k}\ins\N^n$ such that $\pi_i(\nut^{(i),j,k})\eq\nu^{(i),j}$ and
\htt{1.4.5}{}
$$\aligned\bl[\,\nut^{(i),j,k}_i-\nut^{(i),j}_i\,\br]&=k\,\in\,[0,r{-}1],\\\bl\{\nut^{(i),j,k}_i-\nut^{(i),j}_i\br\}&=\be'_{i,j}\defs\bl\{-\nut^{(i),j}_i\br\},\endaligned
\leqno(1.4.5)$$
\par\nin so that
\htt{1.4.6}{}
$$\ell\bl(\nut^{(i),j,k})\eq\al_{i,j}\pl(\be'_{i,j}\pl k)/r.
\leqno(1.4.6)$$
\par\nin Indeed, the coefficient of $\nu_i$ in the linear function $\ell$ coincides with $\tfrac{1}{r}$, since $r\1\ee_i\ins\ell^{-1}(1)$.
Here we may restrict to the case $\al_{i,j}\less 1$ by Remark\,\,\hl{R1.3}{1.3}. (Note that the relation between $[x^{(i),j}]\ins\C\{x_{(i)}\}/(\dd f_{(i)})$ and the eigenvector for $(\al_{i,j},\be_{i,j})$ is {\it not\1} questioned here, since we use only (\hl{1.3.2}{1.3.2}) with $r,r'$ sufficiently large two distinct prime numbers in order to deduce (\hl{1.4.3}{1.4.3}).) This finishes the proof of Proposition\,\,\hl{P1.4}{1.4}.
\par\htt{R1.4a}{}\msn
{\bf Remark\,\,1.4a.} There may be compact faces $\si_k\sst\Gp(f)$ ($k\in[1,m]$) such that the union of $\pi_i(\si_k)$ is a face of $\Gp(f_{(i)})$ with $m\gess2$ (for instance, if $f$ contains $x^2y^2z^a\pl x^4z^b\pl y^4z^b$ with $a\slt b$), see also Remark\,\,\hl{R3.4}{3.4} below.
This is called a {\it fusion of compact faces\1} under a projection. In this case, the argument for the case $\al_{i,j}\sgt1$ would have some difficulty without using Remark\,\,\hl{R1.3}{1.3}.
\par\htt{R1.4b}{}\msn
{\bf Remark\,\,1.4b.} For $n\eq4$, Question~\hl{Q1.4}{1.4} (restricted to the $N$-primitive part) could be solved if the $f_{(i)}$ are {\it simplicial,} the $\SpDef_{f_{(i)}}(t)$ {\it vanish\1} for any $i\ins I_f$, and moreover a fusion of compact faces under a projection as in Remark\,\,\hl{R1.4a}{1.4a} does not occur. We can restrict Question~\hl{Q1.4}{1.4} to the $N$-primitive part in order to get a formula for the $\be_{i,j}$, since the action of $N$ commutes with the local system monodromy.
\ms
We note here the following.
\par\htt{L1.4}{}\msn\vbox{\nin
{\bf Lemma\,\,1.4.} {\it Assume $n\eq 3$. Under the assumption~{\rm (\hl{SC}{SC})}, $f^{-1}(0)$ is reduced, its singular locus is contained in the union of $x_i$-axes for $i\ins I_f$, and the Milnor number $\mu_{(i)}$ of $f_{(i)}$ coincides with $\nu_{(i)}$ which is defined to be the Milnor number of a non-degenerate function $g_i$ such that $\Gp(g_i)\eq\Gp(f_{(i)})$.
Moreover the spectral numbers $\al_{i,j}$ can be determined from the Newton polytope of $f_{(i)}$ using \cite{ste}.}}
\msn
{\it Proof.} Assume there is a factorization $f\eq g_1^ag_2$ in $\C\{x\}$ such that $a$ is an integer greater than 1, $g_1$ is reduced and irreducible, and $g_1^{-1}(0)$ is not contained in $g_2^{-1}(0)$ (where $g_2$ may be an invertible function). Then $(f{+}\ell\1^{ak})^{-1}(0)$ is locally reducible (hence singular) along
$$g_1^{-1}(0)\cap\ell^{-1}(0)\stm g_2^{-1}(0),$$
\par\nin since $a\sgt1$. Here $\ell$ is sufficiently general so that
$$\dim g_1^{-1}(0)\cap\ell^{-1}(0)\cap g_2^{-1}(0)\eq0.$$
\par\nin Hence the above set is non-empty. This contradicts the above assertion that $f{+}\ell\1^r$ is non-degenerate for $r\gg 0$, which implies that it has an isolated singularity at 0, since it is convenient. Thus $f$ must be reduced, hence the singular locus of $f$ has at most dimension 1.
\sk
We now apply Yomdin's formula saying that the difference of the Milnor numbers of $f{+}\ell\1^r$ and $f$ is equal to the sum of the Milnor numbers of the $f_{(i)}$ in \hl{1.3}{1.3} multiplied by $r$ (which follows from (\hl{1.3.1}{1.3.1}) by putting $t\eq1$). Here we have to use the original version as in \cite{ste}, \cite{St2}, where we assume that the singular locus of $f^{-1}(0)$ is a curve (not necessarily a union of lines). In this case the intersection multiplicity $m_i$ of $\ell^{-1}(0)$ with each irreducible component $C_i$ of the curve must be added as a multiplicative factor on the right-hand side of the formula. We thus get
$$\mu_{f+\ell\1^r}\mi\mu_{f+\ell\1^{r'}}=\msum_{i\in I}\,m_i\1\mu_{(i)}(r{-}r')\q\h{if}\,\,\,\,r\sgt r'\gg 0.$$
\par\nin Here $I$ contains $I_f$, and $m_i\eq1$ if $i\ins I_f$.
\sk
On the other hand, we can calculate the difference of Milnor numbers employing \cite{Ko} (together with \cite[Lemma\,\,A.2]{JKSY}), and get
$$\mu_{f+\ell\1^r}\mi\mu_{f+\ell\1^{r'}}=\msum_{i\in I_f}\,\nu_{(i)}(r{-}r')\q\h{if}\,\,\,\,r\sgt r'\gg 0,$$
\par\nin where $\nu_{(i)}$ is as in Lemma\,\,\hl{L1.4}{1.4}. (Note that the determinant of a matrix consisting of three vectors $v_1,v_2,v_3\ins\R^3$ coincides with that for $v_1,v_2-av_1,v_3-bv_1$ for any $a,b\ins\R$. We apply this to the case where $v_1$ is contained in the $x_i$-axis, and $v_2-av_1,v_3-bv_1\ins\{x_i\eq0\}$.)
\sk
We have $\mu_{(i)}\gess\nu_{(i)}$. Indeed, for a deformation $f_{(i)}\pl u\1g_i$ ($u\ins\De_{\ep}$) with $\De_{\ep}$ a disk of radius $\ep\sgt0$, we have $\mu_{f_{(i)}+u\1g_i}\eq\nu_{(i)}$ for $u\ins\De^*_{\ep}$ replacing $\ep$ if necessary.
Comparing the above two equalities, we thus get
$$I\eq I_f,\q\mu_{(i)}\eq\nu_{(i)}.$$
\par\nin The last assertion of Lemma\,\,\hl{L1.4}{1.4} follows from the invariance of spectrum under a $\mu$-constant deformation, see for instance \cite{Va2}.
This finishes the proof of Lemma\,\,\hl{L1.4}{1.4}.
\par\htt{R1.4c}{}\msn
{\bf Remark\,\,1.4c.} It may be possible to prove Lemma\,\,\hl{L1.4}{1.4} for $n\sgt3$ by applying a smooth subdivision of the dual fan of $\Gp(f)$, see \cite[9.1]{Va1} for dual fan. (Here a smooth subdivision means that each cone $\eta$ is a simplex and moreover $\eta\cap\Z^n$ is freely generated by primitive elements in its 1-dimensional subcones.)
\par\htt{1.5}{}\msn
{\bf 1.5.~Monodromy zeta functions.} The first proof of Theorem\,\,\hl{T1}{1} used Varchenko's formula for monodromy zeta functions \cite{Va1}. We review the latter for the convenience of the reader. Assume $f$ is non-degenerate and also simplicial (for simplicity). For an $(n{-}1)$-dimensional compact face $\si$ of $\Gp(f)\sst\R^n$, let $\det(\si)$ be the determinant of the matrix consisting of the coordinates of the vertices of $\si$. Let $\de_{\si}$ be the order of the quotient group of $\Z^n$ by the subgroup generated by the intersection of $\Z^n$ with the affine subspace spanned by $\si$. (If $\si$ is contained in the hyperplane defined by $\msum_{i=1}^n\,a_i\nu_i\eq a_0$ with $a_i\ins\N$, ${\rm GCD}(a_0,\dots,a_n)\eq1$, then $\de_{\si}\eq a_0$.) Set $\mu(\si)\defs|\det(\si)|/\de_{\si}$, and
$$s_{\si}(t)\defs\mu(\si)\,\msum_{k=0}^{\de_{\si}-1}\,t^{k/\de_{\si}},\q M^*_f(t)\defs\msum_{\si\in\CFf^{n-1}}\,s_{\si}(t).$$
\par\nin For a subset $I\sst\{1,\dots,n\}$, set $f^I\defs f|_{\C^I}$, where $\C^I\sst\C^n$ is the subspace defined by $x_i\eq0$ ($i\,{\notin}\,I$). We can define $M^*_{f^I}(t)$ similarly replacing $f$ with $f^I$ and $m$ with $|I|$. Set
$$M_f(t):=\msum_I\,(-1)^{n-|I|}M^*_{f^I}(t),$$
\par\nin where the summation is taken over the subsets $I\sst\{1,\dots,n\}$ including $\emptyset$ and $\{1,\dots,n\}$.
\sk
Define $\phi\,{:}\,\Z[t^{1/e}]\tos\Z[t^{1/e}]$ by $\phi(t^{\al})\defs t^{\al-[\al]}$. The main theorem of \cite{Va1} says that
\htt{1.5.1}{}
$$\phi\bl(\Sp_f(t)\br)\eq M_f(t).
\leqno(1.5.1)$$
\par\nin \par\htt{R1.5}{}\msn
{\bf Remark\,\,1.5.} In the simplicial case we have
\htt{1.5.2}{}
$$s_{\si}(t)=\msum_{\tau\les\si}\,\phi\bl(q_{\si}(t)\br).
\leqno(1.5.2)$$
\par\nin Here the summation is taken over the (not necessarily interior) faces $\tau$ of $\si$ including $\si$ and $\emptyset$ (with $q_{\emptyset}(t)\eq1$), and $\tau\less\si$ means that $\tau$ is a face of $\si$. Indeed, taking an appropriate coordinate change, this can be reduced to the following:
\par\htt{L1.5a}\msn
{\bf Lemma\,\,1.5a.} {\it Assume there is a free $\Z$-submodule $A\subset\Q^n$ of rank $n$ containing $\Z^n$. Then the number of $a=(a_1,\dots,a_n)\in A\cap[0,1)^n$ satisfying the condition $|a|=b/e$ mod $\Z$, is independent of $b\in\Z$, where $|a|:=\msum_{i=1}^n\,a_i$, and $e$ is the positive integer such that $1/e$ is a generator of the subgroup of $\Q$ consisting of $|a|\in\Q$ for $a\in A$.}
\ms
Changing the two conditions on $a$ by $a\in A\cap([0,1)^{n-1}{\times}\R)$ and $|a|=b/e$ (without taking mod $\Z$), and replacing $A$ with the image of $A\cap\{|a|=0\}$ by the projection $\Q^n\to\Q^{n-1}$ forgetting the last coordinate, the above assertion can be reduced to the following (with $n$ replaced by $n{-}1$), since the subgroup $A\cap\{|a|{=}0\}\subset A$ acts on $A\cap\{|a|{=}b/e\}$ freely and transitively so that $A\cap\{|a|{=}b/e\}=A\cap\{|a|{=}0\}+c$ for any $c\in A\cap\{|a|{=}b/e\}$.
\par\htt{L1.5b}{}\msn
{\bf Lemma\,\,1.5b.} {\it If there is a free $\Z$-submodule $A\subset\Q^n$ of rank $n$ containing $\Z^n$, then the number of points of $(A+c)\cap[0,1)^n$ is independent of $c=(c_1,\dots,c_n)\in\R^n$.}
\ms
For the proof of this, we may fix the $c_i$ for $i\ne n$, and move only $c_n$. Then it can be reduced to the bijection
$$(A+c)\cap([0,1)^{n-1}{\times}\{0\})\simeq(A+c)\cap([0,1)^{n-1}{\times}\{1\}).$$
\par\nin \bs\bs\htt{S2}{}
\vbox{\centerline{\bf 2. Proof of Steenbrink's formula}
\bsn
In this section we present a precise proof of Steenbrink's formula for the spectrum in the convenient non-degenerate case.}
\par\htt{2.1}{}\msn
{\bf 2.1.~Steenbrink formula for spectrum.} Let $f$ be a holomorphic function of $n$ variables having a convenient non-degenerate Newton boundary (hence $f$ has an isolated singularity at 0). Here we do {\it not\1} assume $f$ is simplicial. In the notation of \hl{1.2}{1.2}, set
$$A_{\Ga}:=\Gr_{V'_N}^{\ssb}\C\{x\}.$$
\par\nin We denote by $\Gf$ the union of compact faces of $\Gp(f)$. For $\emptyset\less\si\less\Gf$, let $A_{\si}$ be the graded $\C$-subalgebra of $A_{\Ga}$ generated by $[x^{\nu}]$ with $\nu$ contained in the cone $C(\si)$ of $\si$ (where $A_{\emptyset}\eq\C$). We denote by $p_{\si}(t)$ the Hilbert-Poincar\'e series of the graded vector space $A_{\si}$. This is a fractional power series. Define
$$\aligned\qh_{\si}(t)&:=(1-t)^{d(\si)}p_{\si}(t)\q\q\h{with}\q\q d(\si):=\dim C(\si),\\ q_{\si}(t)&:=\msum_{\emptyset\les\tau\les\si}\,(-1)^{d(\si)-d(\tau)}\qh_{\tau}(t).\endaligned$$
\par\nin where the summation is taken over any faces $\tau$ of $\si$ including $\si$ and $\emptyset$. Note that $C(\emptyset)\eq\{0\}$, $\qh_{\emptyset}(t)\eq q_{\emptyset}(t)\eq1$, $d(\si)\eq d_{\si}{+}1$ (even if $\si\eq\emptyset$), and $\qh_{\si}(t)$ is a fractional power polynomial, see \cite{Ko}, \cite{St1}. It is well-known that
\htt{2.1.1}{}
$$\qh_{\si}(t)\eq\msum_{\emptyset\les\tau\les\si}\,q_{\tau}(t).
\leqno(2.1.1)$$
\par\nin (Indeed, we can show that $\msum_{\tau\les\tau'\les\si}\,(-1)^{d_{\si}-d_{\tau'}}\eq0$ for any $\emptyset\less\tau\less\si$.) This definition of $q_{\si}(t)$ is compatible with the one in the introduction if $f$ is simplicial. Set
\htt{2.1.2}{}
$$k(\si)\eq\min\{|I|\mid C(\si)\sst\R^I\},
\leqno(2.1.2)$$
\par\nin where $\R^I\defs\mcap_{i\notin I}\,\{x_i\eq0\}$ for $I\sst\{1,\dots,n\}$. We are interested in the graded $A_{\Ga}$-modules
$$\aligned B_{\Ga}&:=\Gr_{V_N}^{\ssb}\bl(\C\{x\}/(\dd f)\br),\\ B'_{\Ga}&:=\Gr_{V'_N}^{\ssb}\bl(\C\{x\}/(x_1f_1,\dots,x_nf_n)\br),\endaligned$$
\par\nin with $f_i\defs\dd_{x_i}f$. By \cite{exp}, we can calculate the spectrum using (\hl{2.1.3}{2.1.3}) below (where it is {\it not\1} necessary to assume $f$ is \h{\it simplicial\1}), see \cite[Theorem 5.7]{St1}.
\par\htt{T2.1}{}\msn
{\bf Theorem\,\,2.1.} {\it The Poincar\'e polynomial of $B_{\Ga}$ can be described as
\htt{2.1.3}{}
$$p_{B_{\Ga}}(t)=\msum_{\emptyset\les\si\les\Gf}\,(-1)^{n-d(\si)}\,(1{-}t)^{k(\si)-d(\si)}\1\qh_{\si}(t),
\leqno(2.1.3)$$
\par\nin where the summation is taken over any faces of $\Gf$ including $\emptyset$.}
\ms
Note that $V_N$ (but not $V'_N$) is used in the definition of $B_{\Ga}$.
\par\htt{R2.1a}{}\msn
{\bf Remark\,\,2.1a.} Concerning the proof of Theorem\,\,\hl{T2.1}{2.1}, it is written in \cite[p.\,559]{St1} as follows:
\sk
``This can be computed by calculating the Poincar\'e polynomials of $A/(F_0,\cdots,F_n)$ and of all its quotients corresponding to intersections of coordinate hyperplanes in $\R_+^{n+1}.$"
\sk
Here $\,A\eq\Gr_{V'_N}^{\ssb}\C\{x\}$, $\,F_i\eq\Gr_{V'_N}x_if_i\,$ with our notation, and ``This" seems to mean the isomorphism
$$(x_0\cdots x_n)A/(F_0,\dots,F_n)A\cong\Gr\,\C[x_0,\dots,x_n]/(f_0,\dots,f_n),$$
\par\nin assuming that $f$ is a polynomial. However, there is no inclusion
$$(F_0,\dots,F_n)A\sst(x_0\cdots x_n)A$$
\par\nin in the convenient case. Moreover, if $(1,\dots,1)$ belongs to the strict interior of the cone of some maximal-dimensional compact face $\si$, then $(x_0\cdots x_n)A$ must be a subset of $A_{\si}$ because of the {\it algebra structure of\1} $A$. So the assertion seems misstated. The correct one may be as follows: There is a {\it canonical\1} morphism
$$\Gr\,\C[x_0,\dots,x_n]/(f_0,\dots,f_n)\tos A/(F_0,\dots,F_n),$$
\par\nin induced by $\C[x_0,\dots,x_n]\ni g\mapsto\Gr(x_1\cdots x_n)g\ins A$,
and it is furthermore {\it injective,} see also (\hl{2.1.11}{2.1.11}) below.
\sk
The last assertion is, however, {\it highly non-trivial,} and moreover it is quite unclear why this is enough to show the formula (\hl{2.1.3}{2.1.3}).
One problem is that if we take the splitting of the Newton filtration $V'_N$ on the convergent power series $\C\{x\}$ given by {\it monomials}, then this splitting is {\it never\1} compatible with the action of the $F_i$ (unless the Newton polytope has only one maximal-dimensional compact face; in particular, $f$ is weighted homogeneous).
Note also that the Koszul complex for the restriction of $F_i$ $(i\,{\ne}\,j)$ to the coordinate hyperplane defined by $x_j\eq 0$ is a quotient complex of the Koszul complex for the $F_i$, but not a direct factor as a complex.
\par\htt{R2.1b}{}\msn
{\bf Remark\,\,2.1b.} In order to get a precise proof of (\hl{2.1.3}{2.1.3}), it seems necessary to construct a {\it filtered double complex\1} consisting of the Koszul complexes for (the non-zero members of) the restrictions of $x_if_i$ to intersections of coordinate hyperplanes such that its associated single complex is filtered quasi-isomorphic to the filtered Koszul complex for the $f_i$, where the Newton filtration $V_N$ is defined by using the inclusion into the logarithmic complex which is the filtered Koszul complex for the $F_i$.
\sk
More concretely, let $X\sst\C^n$ be a polydisk. Set $X^I\defs\mcap_{i\notin I}\,\{x_i\eq0\}\sst X$ for any subset $I\sst\{1,\dots,n\}$ (including the case $I\eq\emptyset$, where $X^{\emptyset}\eq\{0\}$). Let $\Om^{\log,p}_{X^I}$ be the logarithmic forms on $X^I$ with logarithmic poles along the union of coordinate hyperplanes of $X^I$. These are free sheaves with free generators given by exterior products of the $\ddd x_i/x_i$ ($i\ins I$). The logarithmic forms have the filtration induced by $V'_N$ (but not $V_N$) using the trivialization given by the above free generators.
\sk
We prove the following.
\par\htt{P2.1}{}\msn\vbox{\nin
{\bf Proposition\,\,2.1.} {\it There is a canonical quasi-isomorphism
\htt{2.1.4}{}
$$\bl((\Om^{\ssb}_{X,0},\ddd f\wedge),V_N\br)[n]\simto{\bf s}(\K_0^{\ssb,\ssb},V'_N).
\leqno(2.1.4)$$
\par\nin with $\K^{\ssb,\ssb}$ the {\it residue double complex of logarithmic differential forms\1} such that
$$\K^{j,\ssb}:=\mopl_{|I|=n-j}\,\Om^{\log,\ssb}_{X^I}[n{-}j]\q(j\ins[0,n]),$$
\par\nin where the first differential is induced by residue, and the second by $\ddd f\wedge$.}}
\ms
(We denote by ${\bf s}(\K_0^{\ssb,\ssb})$ the single complex associated to the double complex $\K_0^{\ssb,\ssb}$. Note that $\K_0^{\ssb,\ssb}$ is identified with a double complex consisting of the $($shifted$\1)$ Koszul complexes for the restrictions of $x_if_i$ $(i\ins I)$ to $X^I$ for $I\sst\{1,\dots,n\}$, where the first differential is defined by using the restrictions by $X^{I\setminus\{i\}}\into X^I$ for any $i\ins I$.)
\msn
{\bf Proof of Proposition\,\,\hl{P2.1}{2.1}.} Set $\ddd x_J/x_J\defs\bigwedge_{i\ins J}\ddd x_i/x_i$ (where the exterior product is defined using the natural order of $J\sst\{1,\dots,n\}$). By definition, we have the isomorphisms
$$\Om^{\log,p}_{X^I,0}\eq\mopl_{J\subset I,\,|J|=p}\,\C\{x_I\}\1\ddd x_J/x_J,$$
\par\nin with $\C\{x_I\}\eq\OO_{X^I,0}$.
The first differential is defined by
$$g\1\ddd x_J/x_J\mapsto\msum_{i\in J}\,\ep_{\!J}^{(i)}(g|_{X^{I^{(i)}}})\ddd x_{J^{(i)}}/x_{J^{(i)}}.$$
\par\nin Here $I^{(i)}\defs I\stm\{i\}$, and $\ep_{\!J}^{(i)}\eq(-1)^{k+1}$ if $J\eq\{j_1,\dots,j_p\}$ with $p\eq|J|$, $j_l\slt j_{l+1}$ ($l\ins[1,p-1]$), and $i\eq j_k$.
The second differential is given up to sign by $\ddd f_{(I)}\wedge$, where $f_{(I)}\defs f|_{X^I}$, and the sign is $(-1)^{n-|I|}$. This sign is needed in order that the first and second differentials commute (and not anti-commute). Indeed, we can verify that they anti-commute without this sign. (Note that it is the same as the one needed for the passage from a double complex to the associated single complex.)
We thus get the morphism in (\hl{2.1.4}{2.1.4}).
\sk
For the proof of the quasi-isomorphism, note that any $g_{I,J}\ins\C\{x_I\}\1\ddd x_J/x_J$ ($J\sst I$) is uniquely expressed as
$$\msum_{K\sst J}\,g_{I,J,K}(\ddd x_K/x_K)\wedge\ddd x_{J\setminus K}\q(g_{I,J,K}\ins\C\{x_{I\setminus K}\}),$$
\par\nin using a splitting of the weight filtration on the logarithmic forms (see \cite{De}) by monomials. More precisely, it is obtained by the {\it completed\1} tensor product of the decompositions
$$\C\{x_j\}\1\ddd x_j/x_j\eq\C\{x_j\}\1\ddd x_j\oplus\C\ddd x_j/x_j\q(j\ins J),$$
\par\nin if $J\eq I$. In general, we take its scalar extension by $\C\{x_J\}\into\C\{x_I\}$. (Here ``completed" tensor product means that we apply the scalar extension for the inclusion of the tensor product of the $\C\{x_j\}$ ($j\ins J$) over $\C$ into $\C\{x_J\}$.)
\sk
This induces a direct sum decomposition of the {\it residue\1} complex $\K_0^{\ssb,-q}$ which is indexed by the subsets $K\sst\{1,\dots,n\}$ (including the empty subset) for each $q\ins[0,n]$, since we have by definition
$$\aligned\K_0^{j,-q}&=\mopl_{|I|=n-j}\,\Om^{\log,\,n-j-q}_{X^I,0}\\ &=\mopl_{J\subset I,\,|I|=n-j,\,|J|=n-j-q}\,\C\{x_I\}\1\ddd x_J/x_J.\endaligned$$
\par\nin Here the direct factor indexed by $K\ne\emptyset$ is the Koszul complex associated to the identity morphisms, and hence acyclic. We thus get the canonical quasi-isomorphisms
$$\Om_X^{n-q}\simto\K_0^{\ssb,-q}\q(q\ins[0,n]).$$
\par\nin This implies the quasi-isomorphism (\hl{2.1.4}{2.1.4}) forgetting the filtration.
\sk
As for the relation with the filtration, note that $V'_N$ has a splitting using monomials, and this is compatible with the above direct sum decomposition indexed by $K$. We thus get the filtered quasi-isomorphism (\hl{2.1.4}{2.1.4}), where the filtration $V_N$ on $\Om_{X,0}^{n-q}$ is defined so that the above quasi-isomorphism becomes a filtered quasi-isomorphism; in particular, we have the strict filtered injection
\htt{2.1.5}{}
$$(\Om_{X,0}^p,V_N)\into(\Om_{X,0}^{\log,p},V'_N)\q(p\ins[0,n]),
\leqno(2.1.5)$$
\par\nin that is, $V_N$ is the {\it induced filtration\1} by $V'_N$.
(Note that $V_N$ on $\Om_{X,0}^p$ for $p\ne 0,n$ is not induced from one filtration on $\C\{x\}$ using the free generators $\ddd x_J$ with $|J|\eq p$, since it must be {\it modified depending on\1} $J$.)
This finishes the proof of Proposition\,\,\hl{P2.1}{2.1}.
\msn
{\bf Proof of Theorem\,\,\hl{T2.1}{2.1}.} The second differential of $\K^{\ssb,\ssb}$ is given by $\ddd f\wedge$, and each $\K^{j,\ssb}$ is a direct sum of the (shifted) {\it Koszul complexes\1} associated with the restrictions of the $x_if_i$ as is explained above. We can then apply the quasi-isomorphism of complexes of graded $A_{\Ga}$-modules
\htt{2.1.6}{}
$$A_{\Ga}\simto C_X^{\ssb}\q\,\,\,\,\h{with}\q\,\,\,\,C_X^{\1 j}\defs\msum_{d(\si)=n-j}\,A_{\si}\,\,\,\,(j\ins[0,n]),
\leqno(2.1.6)$$
\par\nin see \cite[Proposition 2.6]{Ko}, after passing to the graded quotients of $V'_N$. Here $\si$ runs over compact faces of $\Gp(f)$ {\it not contained in any coordinate hyperplane of\1} $\R^n$. We have a similar argument with $X$ replaced by $X^I$. These imply the {\it acyclicity\1} except for the top degree:
$$H^k(\Gr_{V'_N}\Om^{\log,\ssb}_{X^I,0},\Gr\,\ddd f\wedge)\eq0\q(k\ne|I|),$$
\par\nin for any $I\sst\{1,\dots,n\}$. We thus get the vanishing
\htt{2.1.7}{}
$$H^q\Gr_{V'_N}\K_0^{j,\ssb}\eq0\q(q\ne 0,\,j\ins\Z),
\leqno(2.1.7)$$
\par\nin and setting
$$B_{\Ga}^{\prime j}\defs H^0\Gr_{V'_N}\K_0^{j,\ssb}\q(j\ins\Z),$$
\par\nin the first differential of $\K^{\ssb,\ssb}$ induces the differential of the complex $B_{\Ga}^{\prime\ssb}$.
By (\hl{2.1.4}{2.1.4}) and (\hl{2.1.7}{2.1.7}) we get the canonical quasi-isomorphisms of complexes
$$(\Gr_{V_N}\Om_{X,0}^{\ssb}[n],\Gr\,\ddd f\wedge)\simto\Gr_{V'_N}\K_0^{\ssb,\ssb}\simto B_{\Ga}^{\prime\ssb},$$
\par\nin where $B_{\Ga}^{\prime j}\eq0$ for $j\notin[0,n]$. These imply the vanishing
\htt{2.1.8}{}
$$H^j(\Gr_{V_N}\Om_{X,0}^{\ssb},\Gr\,\ddd f\wedge)\eq0\q(j\ne n),
\leqno(2.1.8)$$
\par\nin together with the canonical quasi-isomorphism of complexes
\htt{2.1.9}{}
$$B_{\Ga}\simto B_{\Ga}^{\prime\ssb}\q\q\h{with}\q\q B_{\Ga}^{\prime \1 0}\eq B'_{\Ga}.
\leqno(2.1.9)$$
\par\nin Indeed, we have the canonical isomorphisms
\htt{2.1.10}{}
$$\aligned H^n(\Gr_{V_N}\Om_{X,0}^{\ssb},\Gr\,\ddd f\wedge)&=B_{\Ga},\\ H^n(\Gr_{V'_N}\Om_{X,0}^{\log,\ssb},\Gr\,\ddd f\wedge)&=B'_{\Ga},\endaligned
\leqno(2.1.10)$$
\par\nin since $(\Om_{X,0}^{\ssb},V_N,\ddd f\wedge)$ and $(\Om_{X,0}^{\log,\ssb},V'_N,\ddd f\wedge)$ are strict filtered complexes by (\hl{2.1.8}{2.1.8}) and (\hl{2.1.7}{2.1.7}). (In general, a filtered complex $(K^{\ssb},G)$ is {\it strict\1} (that is, $H^jG_pK^{\ssb}\tos H^jK^{\ssb}$ is injective for any $j,p$; hence $\Gr_p^G$ commutes with $H^j$), if $H^jK^{\ssb}\eq H^j\Gr^G_pK^{\ssb}\eq 0$ for $j\ne 0$, $p\ins\Z$, where $G$ is assumed exhaustive. This can be shown by using the long exact sequence associated to the short exact sequences
$$0\tos G_pK^{\ssb}\tos K^{\ssb}\tos K^{\ssb}/G_pK^{\ssb}\tos 0.$$
\par\nin Note that $K^{\ssb}$ is the inductive limit of $G_pK^{\ssb}$, and the inductive limit is an exact functor.)
\sk
The equality (\hl{2.1.3}{2.1.3}) now follows from the two quasi-isomorphisms (\hl{2.1.6}{2.1.6}) and (\hl{2.1.9}{2.1.9}). Concerning the coefficient $(1{-}t)^{k(\si)-d(\si)}$ in (\hl{2.1.3}{2.1.3}), this can be explained as follows:
\sk
For a compact face $\si\less\Gf(f)$ with $d({\si})\slt k(\si)\eq n$, we take the linear subspace $V\sst\R^n$ spanned by $\si$, and take free generators of $\Z^n$ compatible with this subspace, that is, they contain free generators of $V\cap\Z^n$. These give a basis of invariant logarithmic vector fields $\xi_i$ ($i\ins[1,n]$) consisting of linear combinations of $x_i\dd_{x_i}$ (where the {\it dual\1} basis for invariant logarithmic forms is first obtained corresponding to the free generators). So the Koszul complex for the $x_i\dd_{x_i}f_{\si}$ ($i\ins[1,n]$) can be expressed also by the Koszul complex for the $\xi_if_{\si}$ ($i\ins[1,n]$), but $\xi_if_{\si}$ vanishes for $i\notin I_{\si}$, where the $\xi_i$ for $i\ins I_{\si}\sst\{1,\dots,n\}$ correspond to a basis of $V\cap\Z^n$ (with $|I_{\si}|\eq d(\si)$).
\sk
We have a similar argument when $k(\si)\slt n$ (replacing $n$ with $k(\si)$). This produces the coefficient $(1{-}t)^{k(\si)-d(\si)}$ in (\hl{2.1.3}{2.1.3}). As for the sign $(-1)^{n-d(\si)}$, this comes from (\hl{2.1.6}{2.1.6}) and (\hl{2.1.9}{2.1.9}) which give $(-1)^{k(\si)-d(\si)}$ and $(-1)^{n-k(\si)}$ respectively.
Theorem\,\,\hl{T2.1}{2.1} thus follows.
\par\htt{R2.1c}{}\msn
{\bf Remark\,\,2.1c.} By (\hl{2.1.9}{2.1.9}) we get the {\it injectivity\1} of the canonical morphism of graded $A_{\Ga}$-modules
\htt{2.1.11}{}
$$B_{\Ga}\into B'_{\Ga},
\leqno(2.1.11)$$
\par\nin which is induced by (\hl{2.1.5}{2.1.5}) for $p\eq n$, that is, by the multiplication by $x_1\cdots x_n$. (This proves the assertion in the proof of \cite[Theorem 5.7]{St1} which is explained after Theorem\,\,\hl{T2.1}{2.1}.)
\sk
As a corollary of (\hl{2.1.9}{2.1.9}--\hl{2.1.10}{10}), we also see that the graded $A_{\Ga}$-modules $B_{\Ga}$, $B'_{\Ga}$ depend only on the {\it Newton boundary\1} $f_{\Gf}$, where
\htt{2.1.12}{}
$$f_{\Gf}\defs\msum_{\nu\in\N^n\cap\Gf}\1a_{\nu}x^{\nu}\,\,\,\,\h{if}\,\,\,\,f\eq\msum_{\nu\in\N^n}\1a_{\nu}x^{\nu}.
\leqno(2.1.12)$$
\par\nin Note that the graded algebra $A_{\Ga}$ depends only on $\Gf$ (and $\Gp(f)\eq\Gf\pl\R_{\ges 0}^n$).
\par\htt{R2.1d}{}\msn
{\bf Remark\,\,2.1d.} The quasi-isomorphism (\hl{2.1.4}{2.1.4}) can be constructed also by the following argument.
Set $X^{I*}\defs X^I\stm\mcup_{i\in I}\{x_i\eq0\}$ with $j^I:X^{I*}\into X$ the inclusion. There is a canonical quasi-isomorphism of mixed Hodge modules on $X$
$$\Nc^{\ssb}\simto\Q_{h,X}[n]\q\h{with}\q\Nc^{-j}\eq\mopl_{|I|=n-j}\,j^I_!\Q_{h,X^{I*}}[n{-}j],$$
\par\nin see for instance \cite[Remark\,\,2.4\,(iv)]{FFS}. Taking the dual, this induces a the canonical quasi-isomorphism of mixed Hodge modules
$$\Q_{h,X}[n]\simto\M^{\ssb}\q\h{with}\q\M^j\eq\mopl_{|I|=n-j}\,j^I_*\Q_{h,X^{I*}}(-j)[n{-}j].$$
\par\nin Applying the filtered de Rham functor, we then get a canonical filtered quasi-isomorphism
$$(\Om_X^{\ssb}[n],F)\simto(\K^{\ssb,\ssb},F),$$
\par\nin where the Hodge filtration $F$ on $\Om_X^{\ssb}[n]$ is defined by the truncation $\si^{\ges p}$ (see \cite{De}), and similarly for each $\K^{j,\ssb}\eq\mopl_{|I|=n-j}\,\Om^{\log,\ssb}_{X^I}[n{-}j]$ ($j\ins[0,n]$). Note that the second differential of $\K^{\ssb,\ssb}$ is given by $\ddd$ (instead of $\ddd f\wedge$), and the first differential is induced from the residue morphisms (up to sign).
\sk
Taking $\Gr_F^{-p}$ for $p\ins[0,n]$, it induces the quasi-isomorphisms
$$\Om_X^{n-p}\simto\K^{\ssb,-p}\q\h{with}\q \K^{j,-p}\eq\mopl_{|I|=n-j}\,\Om_{X^I}^{\log,\,n-j-p},$$
\par\nin since $F^{-p}$ is defined by $\si^{\ges-p}$. Here the differential on the right-hand side is induced by the residue morphisms (up to sign). We can verify that this is compatible with $\ddd f\wedge$. So we get the canonical quasi-isomorphism (\hl{2.1.4}{2.1.4}).
\par\htt{R2.1e}{}\msn
{\bf Remark\,\,2.1e.}
If one tries to prove (\hl{2.1.3}{2.1.3}) using the weight filtration $W$ on the logarithmic complex $\Omega_X^{\log,\ssb}$, then there is a problem of strict bifiltered complexes for $(\Omega_X^{\log,\ssb};V,W)$, that is, the commutation of the three operations $\Gr^W_{\ssb}$, $\Gr_V^{\ssb}$, $H^{\ssb}$, which is essentially the problem of three filtrations. This argument seems to be more complicated than the one using the residue double complex, since we would have to prove also the compatibility of the filtration $V$ on the usual differential forms $\Omega_{X^I}^p$, which can be induced from various $\Omega_{X^J}^{\log,\ssb}$ via $\Gr^W_{\ssb}$ for any $J\supset I$.
\par\htt{R2.1f}{}\msn
{\bf Remark\,\,2.1f.} It is not necessarily easy to prove the coincidence of the $V$-filtration on the {\it Brieskorn lattice\1} $H''_f$ with the filtration induced by $V_N$ on $\C\{x\}$ in (\hl{1.2.4}{1.2.4}). This does not immediately follow from the assertion that the filtration $V_N$ on the {\it microlocal Gauss-Manin complex\1} induces the $V$-filtration on the {\it Gauss-Manin system.} (This is shown in \cite{exp}.) We need also the {\it bistrictness\1} of the microlocal Gauss-Manin complex filtered by $F,V_N$, which is also shown in \cite{exp} essentially and implies that the three functors $H^0$, $F_p/F_q$, $V^{\al}/V^{\be}$ commute with each other for $p\sgt q$, $\al\slt\be$, including the case $q\eq{-\infty}$, $\be\eq{+\infty}$. Indeed, using the coincidence on the Gauss-Manin system together with the bistrictness of the microlocal Gauss-Manin complex, the assertion is reduced to that $\dd_t^{-1}\omega\ins V_N^{\al+1}H''_f$ if $\omega\ins V_N^{\al}H''_f$ (by the definition of the microlocal Hodge filtration). Here we may assume that $\omega$ is represented by $g\ddd x$ with $g$ a {\it monomial\1} by the definition of $V_N$. Then $\dd_t^{-1}\omega$ is represented by $x_i(\dd_{x_i}f)g\1\ddd x$ {\it up to constant multiple,} and the assertion follows from the inclusion $\Gp(f)\pl\be\Gp(f)\sst(\be{+}1)\Gp(f)$ for $\be\ins\R_{>0}$, see \cite{Ko}. (A similar argument applies to the semi-weighted-homogeneous case, see also \cite[Remark 2.2d]{exa}.)
\par\htt{2.2}{}\msn
{\bf 2.2.~Determination of Newton polyhedra for $n\eq 3$.} It is not necessarily easy to {\it determine\1} the Newton polyhedron {\it even\1} in the 3 variable case. Indeed, if there are 4 points $p_i\ins\R_{\ges 0}^3$ ($i\ins[1,4]$) such that $[p_i,p_j]$ for $i\ins\{1,2\}$, $j\ins\{3,4\}$ are faces of the Newton polytope, it is not necessarily trivial to decide whether $[p_1,p_2]$ or $[p_3,p_4]$ is a face of it. Here $[p,q]$ means the convex hull of $\{p,q\}\sst\R^3$. To solve the problem, we may have to use the following:
\par\htt{L2.2}{}\msn
{\bf Lemma\,\,2.2.} {\it Let $p_i$ $(i\ins[1,4])$ be $4$ points on $\R^3$ not lying on one affine plane. Assume any three of them are linearly independent. Let $\ell_i$ be the line joining $p_{2i}$ and $p_{2i-1}$ $(i\eq 1,2)$. Let $\ell_0$ be the intersection of the cones of $\ell_1$ and $\ell_2$. Assume there is an intersection point $u_i\ins\ell_0\cap\ell_i$ $(i\eq 1,2)$, and moreover $u_2\eq ru_1$ with $r\sgt0$. Then we have $r\slt1$ if and only if}
\htt{2.2.1}{}
$$\det(p_1{-}p_3,p_2{-}p_3,p_4{-}p_3)\det(p_1,p_2,p_4{-}p_3)\sgt0.
\leqno(2.2.1)$$
\par\nin \msn
{\it Proof.} Changing the coordinates of $\R^3$, we may assume that $p_1,p_2,p_3$ are respectively $\ee_1,\ee_2,\ee_3$ with $\ee_i$ the \h{$i$\1th} unit vector. Write $p_4\eq(a,b,c)\ins\R^3$.
The intersection of the line joining $p_3,p_4$ and the plane containing $0,p_1,p_2$ is given by $\bl(\tfrac{a}{1{-}c},\tfrac{b}{1{-}c},0\br)$. This implies that $r\eq\tfrac{a{+}b}{1{-}c}$. Hence we have $r\slt1$ if and only if $(a{+}b{+}c{-}1)(c{-}1)\sgt0$. The assertion then follows from the equalities
$$\det\Biggl(\!\!\begin{array}{rrc}1&\!\!0&\!\!a\\ 0&\!\!1&\!\!b\\ -1&\!\!-1&\!\!c{-}1\end{array}\!\Biggl)\eq a\pl b\pl c\mi 1,\q\det\Biggl(\!\begin{array}{ccc}1&0&\!a\\ 0&1&\!b\\ 0&0&\!c{-}1\end{array}\!\Biggl)\eq c\mi 1.$$
\par\nin This finishes the proof of Lemma\,\,\hl{L2.2}{2.2}.
\par\htt{R2.2a}{}\msn
{\bf Remark\,\,2.2a.} The determinant of $(p_1{-}p_3,p_2{-}p_3,p_4{-}p_3)$ vanishes if and only if the 4 points are on one affine plane. So only the problem of sign remains.
\par\htt{E2.2}{}\msn
{\bf Example~2.2.} For $j\eq17,19$, set
$$f_j\eq x^j{+}y^{15}{+}z^{15}{+}x^3y^6{+}x^6y^3{+}x^4z^8{+}x^8z^4{+}y^2z^4{+}y^4z^2.$$
\par\nin This example shows that the defective spectral numbers may depend on the {\it combinatorics\1} of the Newton polytope. Indeed, calculating the $\ga_{\si}$, we can get by Theorem\,\,\hl{T1}{1} that
$$\SpDef_{f_j}(t)=\begin{cases}t\1\bl(4\1t\1^{1/2}+\msum_{i=1}^{2}\,t\1^{i/3}+\msum_{i=1}^{16}\,t\1^{i/17}\1\br)&\h{if}\,\,\,\,j\eq17,\\ t\1\bl(3\1t\1^{1/2}+2\1\msum_{i=1}^{2}\,t\1^{i/3}+\msum_{i=1}^{3}\,t\1^{i/4}\1\br)&\h{if}\,\,\,\,j\eq19.\end{cases}$$
\par\nin Here the $q_{\si}(t)$ for $\si\eq(0,4,2)$, $(0,2,4)$, $(6,3,0)$, $(8,0,4)$, $(17,0,0)$ are given respectively by
$$t\1^{1/2},\q t\1^{1/2},\q\msum_{i=1}^{2}\,t\1^{i/3},\q\msum_{i=1}^{3}\,t\1^{i/4},\q\msum_{i=1}^{16}\,t\1^{i/17}.$$
\par\nin (Note that $f_{18}$ is not simplicial.) These agree with calculations using a computer mentioned in Remark\,\,\hl{R3.1a}{3.1a} below. Here we can count the $\ga_{\si}$ using the pictures below:
\sk
$$\setlength{\unitlength}{.4cm}
\begin{picture}(13,8)
\put(3,0.9){\line(-1,-1){.9}}
\put(3,6.9){\line(-1,2){.5}}
\put(9,0.9){\line(2,-1){1}}
\put(3,0.9){\line(0,1){6}}
\put(3,0.9){\line(1,0){6}}
\put(9,0.9){\line(-1,1){6}}
\put(3,4.9){\line(1,0){2}}
\put(3,4.9){\line(2,-1){4}}
\put(3,2.9){\line(1,0){4}}
\put(3,2.9){\line(3,-1){6}}
\put(3,2.9){\line(2,-1){4}}
\put(3,2.9){\line(1,-1){2}}
\put(9.3,0.9){$\scriptstyle (17,0,0)$}
\put(7.3,2.9){$\scriptstyle (6,3,0)$}
\put(5.3,4.9){$\scriptstyle (3,6,0)$}
\put(.3,0.9){$\scriptstyle (0,0,15)$}
\put(.5,2.7){$\scriptstyle (0,2,4)$}
\put(.5,4.7){$\scriptstyle (0,4,2)$}
\put(.3,6.5){$\scriptstyle (0,15,0)$}
\put(3.7,.2){$\scriptstyle (4,0,8)$}
\put(6.2,.2){$\scriptstyle (8,0,4)$}
\end{picture}\q
\begin{picture}(13,7.9)
\put(3,0.9){\line(-1,-1){.9}}
\put(3,6.9){\line(-1,2){.5}}
\put(9,0.9){\line(2,-1){1}}
\put(3,0.9){\line(0,1){6}}
\put(3,0.9){\line(1,0){6}}
\put(9,0.9){\line(-1,1){6}}
\put(3,4.9){\line(1,0){2}}
\put(3,4.9){\line(2,-1){4}}
\put(3,2.9){\line(1,0){4}}
\put(7,0.9){\line(0,1){2}}
\put(3,2.9){\line(2,-1){4}}
\put(3,2.9){\line(1,-1){2}}
\put(9.3,0.9){$\scriptstyle (19,0,0)$}
\put(7.3,2.9){$\scriptstyle (6,3,0)$}
\put(5.3,4.9){$\scriptstyle (3,6,0)$}
\put(.3,0.9){$\scriptstyle (0,0,15)$}
\put(.5,2.7){$\scriptstyle (0,2,4)$}
\put(.5,4.7){$\scriptstyle (0,4,2)$}
\put(.3,6.5){$\scriptstyle (0,15,0)$}
\put(3.7,.2){$\scriptstyle (4,0,8)$}
\put(6.2,.2){$\scriptstyle (8,0,4)$}
\end{picture}$$
\par\nin \msn
These two triangulations of the triangle $\Ga\defs H\cap\R_{\ges0}^3$ are given by intersecting the cones $C(\si)$ with the plane $H\defs\bl\{\msum_{i=1}^3\,\nu_i\eq1\br\}\sst\R^3$ for compact faces $\si$ of $\Gp({f_{17}})$ and $\Gp({f_{19}})$ respectively. (Note, however, that these triangulations cannot determine the toric varieties associated with the Newton polytopes, since some information is lost by passing from $\si$ to $C(\si)$.)
Here $\ga_{\si}$ for $\si\ins\CFf^0$ is given by the number of edges containing the point of $\Ga$ corresponding to $\si$ (since $\dim\Ga\eq2$).
Three edges outside $\Ga$ represent the axes of $\R^3$. Notice that the edge joining the points corresponding to $(0,2,4),(17,0,0)\ins\Gp(f_{17})$ is changed to the one for $(6,3,0),(8,0,4)\ins\Gp(f_{19})$. Consequently some $\ga_{\si}$ are also changed.
\par\htt{R2.2b}{}\msn
{\bf Remark\,\,2.2b.} For 4 points $(6,3,0)$, $(8,0,4)$, $(0,2,4)$, $(j,0,0)$ in Example~\hl{E2.2}{2.2} just above, the two determinants in (\hl{2.2.1}{2.2.1}) are respectively given as follows:
$$\det\Biggl(\!\!\begin{array}{rrr}
6&\!\!8&\!\!j\\
1&\!\!-2&\!\!-2\\
-4&\!\!0&\!\!-4\end{array}\!\Biggl)\eq8\1(18\mi j),\q
\det\Biggl(\!\begin{array}{rrr}
6&8&\!\!j\\
3&0&\!\!-2\\
0&4&\!\!-4\end{array}\!\Biggl)\eq12\1(j\pl 12).$$
\par\nin \par\htt{2.3}{}\msn
{\bf 2.3.~Topological classification of Newton boundaries for $|\CFf^0|\eq6$, $n\eq3$.} In the convenient case with $n\eq3$, three points of $\CFf^0$ are contained in the coordinate axes, and the topological classification of simplicial convenient Newton boundaries (as triangulated spaces) is rather easy for $|\CFf^0|\less 5$. In the case $|\CFf^0|\eq6$, we have the following topological types up to permutation of coordinates (see Example\,\,\hl{E2.2}{2.2} for the meaning of pictures).
\ms
\scalebox{0.57}{$\setlength{\unitlength}{.4cm}
\begin{picture}(10,9)
\put(1,1){\line(-1,-1){.9}}
\put(1,7){\line(-1,2){.5}}
\put(7,1){\line(2,-1){1}}
\put(1,1){\line(0,1){6}}
\put(1,1){\line(1,0){6}}
\put(7,1){\line(-1,1){6}}
\put(1,1){\line(1,1){3}}
\put(4,1){\line(0,1){3}}
\put(1,4){\line(1,0){3}}
\end{picture}
\begin{picture}(10,9)
\put(1,1){\line(-1,-1){.9}}
\put(1,7){\line(-1,2){.5}}
\put(7,1){\line(2,-1){1}}
\put(1,1){\line(0,1){6}}
\put(1,1){\line(1,0){6}}
\put(7,1){\line(-1,1){6}}
\put(4,1){\line(-1,1){3}}
\put(4,1){\line(0,1){3}}
\put(1,4){\line(1,0){3}}
\end{picture}
\begin{picture}(10,9)
\put(1,1){\line(-1,-1){.9}}
\put(1,7){\line(-1,2){.5}}
\put(7,1){\line(2,-1){1}}
\put(1,1){\line(0,1){6}}
\put(1,1){\line(1,0){6}}
\put(7,1){\line(-1,1){6}}
\put(1,1){\line(2,1){4}}
\put(3,1){\line(1,1){2}}
\put(1,1){\line(1,2){2}}
\end{picture}
\begin{picture}(10,9)
\put(1,1){\line(-1,-1){.9}}
\put(1,7){\line(-1,2){.5}}
\put(7,1){\line(2,-1){1}}
\put(1,1){\line(0,1){6}}
\put(1,1){\line(1,0){6}}
\put(7,1){\line(-1,1){6}}
\put(3,1){\line(0,1){4}}
\put(3,1){\line(1,1){2}}
\put(1,1){\line(1,2){2}}
\end{picture}
\begin{picture}(10,9)
\put(1,1){\line(-1,-1){.9}}
\put(1,7){\line(-1,2){.5}}
\put(7,1){\line(2,-1){1}}
\put(1,1){\line(0,1){6}}
\put(1,1){\line(1,0){6}}
\put(7,1){\line(-1,1){6}}
\put(3,1){\line(0,1){4}}
\put(3,1){\line(1,1){2}}
\put(3,1){\line(-1,3){2}}
\end{picture}
\begin{picture}(10,9)
\put(1,1){\line(-1,-1){.9}}
\put(1,7){\line(-1,2){.5}}
\put(7,1){\line(2,-1){1}}
\put(1,1){\line(0,1){6}}
\put(1,1){\line(1,0){6}}
\put(7,1){\line(-1,1){6}}
\put(1,1){\line(1,1){3}}
\put(1,1){\line(2,1){4}}
\put(1,1){\line(1,2){2}}
\end{picture}
\begin{picture}(10,9)
\put(1,1){\line(-1,-1){.9}}
\put(1,7){\line(-1,2){.5}}
\put(7,1){\line(2,-1){1}}
\put(1,1){\line(0,1){6}}
\put(1,1){\line(1,0){6}}
\put(7,1){\line(-1,1){6}}
\put(1,1){\line(1,1){1.5}}
\put(2.5,2.5){\line(0,1){3}}
\put(2.5,2.5){\line(1,0){3}}
\put(2.5,2.5){\line(-1,3){1.5}}
\put(2.5,2.5){\line(3,-1){4.5}}
\end{picture}$}
\ms
\scalebox{0.57}{$\setlength{\unitlength}{.4cm}
\begin{picture}(10,9)
\put(1,1){\line(-1,-1){.9}}
\put(1,7){\line(-1,2){.5}}
\put(7,1){\line(2,-1){1}}
\put(1,1){\line(0,1){6}}
\put(1,1){\line(1,0){6}}
\put(7,1){\line(-1,1){6}}
\put(1,1){\line(1,1){1.5}}
\put(2.5,2.5){\line(0,1){3}}
\put(2.5,2.5){\line(1,0){3}}
\put(1,1){\line(1,3){1.5}}
\put(2.5,2.5){\line(3,-1){4.5}}
\end{picture}
\begin{picture}(10,9)
\put(1,1){\line(-1,-1){.9}}
\put(1,7){\line(-1,2){.5}}
\put(7,1){\line(2,-1){1}}
\put(1,1){\line(0,1){6}}
\put(1,1){\line(1,0){6}}
\put(7,1){\line(-1,1){6}}
\put(1,1){\line(1,1){1.5}}
\put(2.5,2.5){\line(0,1){3}}
\put(2.5,2.5){\line(1,0){3}}
\put(1,1){\line(3,1){4.5}}
\put(1,1){\line(1,3){1.5}}
\end{picture}
\begin{picture}(10,9)
\put(1,1){\line(-1,-1){.9}}
\put(1,7){\line(-1,2){.5}}
\put(7,1){\line(2,-1){1}}
\put(1,1){\line(0,1){6}}
\put(1,1){\line(1,0){6}}
\put(7,1){\line(-1,1){6}}
\put(1,1){\line(1,3){1.5}}
\put(1,1){\line(1,1){3}}
\put(4,2){\line(0,1){2}}
\put(1,1){\line(3,1){3}}
\put(4,2){\line(3,-1){3}}
\end{picture}
\begin{picture}(10,9)
\put(1,1){\line(-1,-1){.9}}
\put(1,7){\line(-1,2){.5}}
\put(7,1){\line(2,-1){1}}
\put(1,1){\line(0,1){6}}
\put(1,1){\line(1,0){6}}
\put(7,1){\line(-1,1){6}}
\put(3,3){\line(-1,0){2}}
\put(1,3){\line(3,-1){6}}
\put(3,3){\line(-1,2){2}}
\put(3,3){\line(2,-1){4}}
\put(1,3){\line(1,-1){2}}
\end{picture}
\begin{picture}(10,9)
\put(1,1){\line(-1,-1){.9}}
\put(1,7){\line(-1,2){.5}}
\put(7,1){\line(2,-1){1}}
\put(1,1){\line(0,1){6}}
\put(1,1){\line(1,0){6}}
\put(7,1){\line(-1,1){6}}
\put(3,3){\line(-1,0){2}}
\put(3,3){\line(0,-1){2}}
\put(3,3){\line(-1,2){2}}
\put(3,3){\line(2,-1){4}}
\put(1,3){\line(1,-1){2}}
\end{picture}
\begin{picture}(10,9)
\put(1,1){\line(-1,-1){.9}}
\put(1,7){\line(-1,2){.5}}
\put(7,1){\line(2,-1){1}}
\put(1,1){\line(0,1){6}}
\put(1,1){\line(1,0){6}}
\put(7,1){\line(-1,1){6}}
\put(3,3){\line(-1,0){2}}
\put(3,3){\line(0,-1){2}}
\put(3,3){\line(-1,2){2}}
\put(3,3){\line(2,-1){4}}
\put(1,1){\line(1,1){2}}
\end{picture}
\begin{picture}(10,9)
\put(1,1){\line(-1,-1){.9}}
\put(1,7){\line(-1,2){.5}}
\put(7,1){\line(2,-1){1}}
\put(1,1){\line(0,1){6}}
\put(1,1){\line(1,0){6}}
\put(7,1){\line(-1,1){6}}
\put(1,5.5){\line(4,-3){6}}
\put(2.5,2.5){\line(3,-1){4.5}}
\put(2.5,2.5){\line(-1,2){1.5}}
\put(2.5,2.5){\line(0,-1){1.5}}
\put(1,1){\line(1,1){1.5}}
\end{picture}$}
\ms
\scalebox{0.57}{$\setlength{\unitlength}{.4cm}
\begin{picture}(10,9)
\put(1,1){\line(-1,-1){.9}}
\put(1,7){\line(-1,2){.5}}
\put(7,1){\line(2,-1){1}}
\put(1,1){\line(0,1){6}}
\put(1,1){\line(1,0){6}}
\put(7,1){\line(-1,1){6}}
\put(1,5.5){\line(4,-3){6}}
\put(2.5,2.5){\line(3,-1){4.5}}
\put(2.5,2.5){\line(-1,2){1.5}}
\put(2.5,2.5){\line(0,-1){1.5}}
\put(1,5.5){\line(1,-3){1.5}}
\end{picture}
\begin{picture}(10,9)
\put(1,1){\line(-1,-1){.9}}
\put(1,7){\line(-1,2){.5}}
\put(7,1){\line(2,-1){1}}
\put(1,1){\line(0,1){6}}
\put(1,1){\line(1,0){6}}
\put(7,1){\line(-1,1){6}}
\put(1,7){\line(1,-4){1}}
\put(1,1){\line(1,1){3}}
\put(1,1){\line(2,1){2}}
\put(1,1){\line(1,2){1}}
\put(2,3){\line(2,1){2}}
\put(3,2){\line(1,2){1}}
\put(7,1){\line(-4,1){4}}
\end{picture}
\begin{picture}(10,9)
\put(1,1){\line(-1,-1){.9}}
\put(1,7){\line(-1,2){.5}}
\put(7,1){\line(2,-1){1}}
\put(1,1){\line(0,1){6}}
\put(1,1){\line(1,0){6}}
\put(7,1){\line(-1,1){6}}
\put(1,7){\line(1,-4){1}}
\put(2,3){\line(1,-1){1}}
\put(1,1){\line(2,1){2}}
\put(1,1){\line(1,2){1}}
\put(2,3){\line(2,1){2}}
\put(3,2){\line(1,2){1}}
\put(7,1){\line(-4,1){4}}
\end{picture}
\begin{picture}(10,9)
\put(1,1){\line(-1,-1){.9}}
\put(1,7){\line(-1,2){.5}}
\put(7,1){\line(2,-1){1}}
\put(1,1){\line(0,1){6}}
\put(1,1){\line(1,0){6}}
\put(7,1){\line(-1,1){6}}
\put(1,7){\line(1,-4){1}}
\put(2,3){\line(2,-1){2}}
\put(1,1){\line(3,1){3}}
\put(1,1){\line(1,2){1}}
\put(1,7){\line(3,-5){3}}
\put(4,2){\line(0,1){2}}
\put(7,1){\line(-3,1){3}}
\end{picture}
\begin{picture}(10,9)
\put(1,1){\line(-1,-1){.9}}
\put(1,7){\line(-1,2){.5}}
\put(7,1){\line(2,-1){1}}
\put(1,1){\line(0,1){6}}
\put(1,1){\line(1,0){6}}
\put(7,1){\line(-1,1){6}}
\put(1,7){\line(1,-5){1}}
\put(2,2){\line(2,0){2}}
\put(2,2){\line(5,-1){5}}
\put(1,1){\line(1,1){1}}
\put(1,7){\line(3,-5){3}}
\put(4,2){\line(0,1){2}}
\put(7,1){\line(-3,1){3}}
\end{picture}
\begin{picture}(10,9)
\put(1,1){\line(-1,-1){.9}}
\put(1,7){\line(-1,2){.5}}
\put(7,1){\line(2,-1){1}}
\put(1,1){\line(0,1){6}}
\put(1,1){\line(1,0){6}}
\put(7,1){\line(-1,1){6}}
\put(1,7){\line(1,-5){1}}
\put(2,2){\line(2,0){2}}
\put(2,2){\line(5,-1){5}}
\put(1,1){\line(1,1){1}}
\put(2,2){\line(1,1){2}}
\put(4,2){\line(0,1){2}}
\put(7,1){\line(-3,1){3}}
\end{picture}
\begin{picture}(10,9)
\put(1,1){\line(-1,-1){.9}}
\put(1,7){\line(-1,2){.5}}
\put(7,1){\line(2,-1){1}}
\put(1,1){\line(0,1){6}}
\put(1,1){\line(1,0){6}}
\put(7,1){\line(-1,1){6}}
\put(1,7){\line(1,-4){1}}
\put(2,3){\line(1,0){3}}
\put(1,1){\line(2,1){4}}
\put(1,1){\line(1,2){1}}
\put(1,1){\line(1,1){2}}
\put(1,7){\line(1,-2){2}}
\end{picture}$}
\ms
\scalebox{0.57}{$\setlength{\unitlength}{.4cm}
\begin{picture}(10,9)
\put(1,1){\line(-1,-1){.9}}
\put(1,7){\line(-1,2){.5}}
\put(7,1){\line(2,-1){1}}
\put(1,1){\line(0,1){6}}
\put(1,1){\line(1,0){6}}
\put(7,1){\line(-1,1){6}}
\put(3,3){\line(1,0){2}}
\put(3,4){\line(0,-1){1}}
\put(3,4){\line(2,-1){2}}
\put(3,4){\line(-2,3){2}}
\put(3,4){\line(-2,-3){2}}
\put(1,1){\line(2,1){4}}
\put(1,1){\line(1,1){2}}
\end{picture}
\begin{picture}(10,9)
\put(1,1){\line(-1,-1){.9}}
\put(1,7){\line(-1,2){.5}}
\put(7,1){\line(2,-1){1}}
\put(1,1){\line(0,1){6}}
\put(1,1){\line(1,0){6}}
\put(7,1){\line(-1,1){6}}
\put(1,7){\line(1,-4){1}}
\put(2,3){\line(1,0){3}}
\put(2,3){\line(1,1){1}}
\put(1,1){\line(2,1){4}}
\put(1,1){\line(1,2){1}}
\put(3,4){\line(2,-1){2}}
\put(1,7){\line(2,-3){2}}
\end{picture}
\begin{picture}(10,9)
\put(1,1){\line(-1,-1){.9}}
\put(1,7){\line(-1,2){.5}}
\put(7,1){\line(2,-1){1}}
\put(1,1){\line(0,1){6}}
\put(1,1){\line(1,0){6}}
\put(7,1){\line(-1,1){6}}
\put(1,1){\line(1,1){2}}
\put(7,1){\line(-5,1){5}}
\put(1,7){\line(1,-5){1}}
\put(7,1){\line(-3,1){4.5}}
\put(1,7){\line(1,-3){1.5}}
\put(7,1){\line(-2,1){4}}
\put(1,7){\line(1,-2){2}}
\end{picture}
\begin{picture}(10,9)
\put(1,1){\line(-1,-1){.9}}
\put(1,7){\line(-1,2){.5}}
\put(7,1){\line(2,-1){1}}
\put(1,1){\line(0,1){6}}
\put(1,1){\line(1,0){6}}
\put(7,1){\line(-1,1){6}}
\put(1,1){\line(1,1){1.5}}
\put(7,1){\line(-2,1){3}}
\put(1,7){\line(1,-2){1.5}}
\put(2.5,2.5){\line(1,0){1.5}}
\put(2.5,2.5){\line(0,1){1.5}}
\put(2.5,4){\line(1,-1){1.5}}
\put(1,1){\line(2,1){3}}
\put(1,1){\line(1,2){1.5}}
\put(2.5,4){\line(3,-2){4.5}}
\end{picture}
\begin{picture}(10,9)
\put(1,1){\line(-1,-1){.9}}
\put(1,7){\line(-1,2){.5}}
\put(7,1){\line(2,-1){1}}
\put(1,1){\line(0,1){6}}
\put(1,1){\line(1,0){6}}
\put(7,1){\line(-1,1){6}}
\put(1,1){\line(1,1){1.5}}
\put(7,1){\line(-2,1){3}}
\put(1,7){\line(1,-2){1.5}}
\put(2.5,2.5){\line(1,0){1.5}}
\put(2.5,2.5){\line(0,1){1.5}}
\put(2.5,4){\line(1,-1){1.5}}
\put(1,1){\line(2,1){3}}
\put(1,7){\line(1,-3){1.5}}
\put(2.5,4){\line(3,-2){4.5}}
\end{picture}
\begin{picture}(10,9)
\put(1,1){\line(-1,-1){.9}}
\put(1,7){\line(-1,2){.5}}
\put(7,1){\line(2,-1){1}}
\put(1,1){\line(0,1){6}}
\put(1,1){\line(1,0){6}}
\put(7,1){\line(-1,1){6}}
\put(1,1){\line(1,1){2}}
\put(1,1){\line(2,1){2}}
\put(1,1){\line(1,2){1}}
\put(2,3){\line(1,0){1}}
\put(3,2){\line(0,1){1}}
\put(7,1){\line(-2,1){4}}
\put(7,1){\line(-4,1){4}}
\put(1,7){\line(1,-2){2}}
\put(1,7){\line(1,-4){1}}
\end{picture}$}
\msn
It does not seem necessarily easy to verify that these are complete. Some of them can be obtained from others by {\it changing triangulations of subquadrilaterals,} see for instance the first and second types (and also the third, fourth, and fifth). It does not seem trivial either to construct a polynomial corresponding to each type (especially when there is a vertex contained in many edges). If a triangulation of $\Ga$ induces a triangulation of a subtriangle $\Ga'$ consisting of three edges joining an interior point with three vertices of $\Ga'$, then the problem can be reduced to the case where the three edges are removed, since one can add a vertex very close to the 2-dimensional compact face corresponding to $\Ga'$ of the simplified Newton polyhedron. (The last type contains two such subtriangles.) For $|\CFfin^0|\eq3$, the topological type is determined by the $\ga_{\si}\mi2$ (which is equal to the number of $\si'\ins\CFfin^0$ connected with $\si$ by an edge) for $\si\ins\CFf^0\stm\CFfin^0$. Note that $\ga_{\si}$ is given by the number of edges containing $\si\ins\CFf^0$ when $n\eq3$.
\par\htt{2.4}{}\msn
{\bf 2.4.~Examples for $n\eq 4$.} Assume more generally $n\gess 3$. Set $\ff_i\defs\ob{-}\,\ee_i$, where $\ob\defs\msum_{i=1}^n\,\ee_i$ and $\ee_i$ is the \h{$i$\1th} unit vector $(i\ins[1,n]$). Let $(a,b,c)$ be the minimal triple of positive integers in lexicographic order satisfying the following condition: There is a polynomial $f\ins\C[x_1,\dots,x_n]$ such that $f$ is stable by the action of the symmetric group ${\mathfrak S}_n$, the vertices of the Newton polytope $\Gp(f)$ contain $\,a\ob$, $b\1\ff_i$, $c\1\ee_i$ ($i\ins[1,n]$), and moreover the convex hull of $\bl\{a\1\ob$, $b\1\ff_i$ ($i\ins[2,n])\br\}$ is a face of $\Gp(f)$. It may be expected in general that
$$(a,b,c)\eq(2n{-}3,2n{-}1,(2n{-}3)(2n{-}1)\pl1).$$
\par\nin Indeed, the hyperplane containing the above convex hull is defined by the equation
$$\bl((n{-}1)a\mi(n{-}2)b\br)\nu_1+(b\mi a)\1\msum_{i=2}^n\,\nu_i\eq ab,$$
\par\nin where $\nu\eq(\nu_i)\sst\R^n$. This implies that
$$\tfrac{n-2}{n-1}\slt\tfrac{a}{b}\slt\tfrac{n-1}{n}\,\,({<}\,1),$$
\par\nin where the second inequality follows by comparing the degrees of monomials corresponding to $a\1\ob$ and $b\1\ff_i$. Taking the differences with 1, these inequalities imply that
$$\tfrac{1}{n-1}\sgt\tfrac{b-a}{b}\sgt\tfrac{1}{n},\q\h{hence}\q n{-}2\slt\tfrac{a}{d}\slt n{-}1,$$
\par\nin where $d\defs b\mi a$ so that $\tfrac{b}{b-a}\eq\tfrac{a}{d}\pl1$. We then get that
$$d\eq 2,\q a\eq2n{-}3,\q b\eq2n{-}1,$$
\par\nin as {\it numerical candidates.} (Here it seems quite non-trivial to determine $\Gp(f)$.) We would have $c\eq ab\pl 1$, since we must get the inequality $c\sgt ab$ by substituting $\nu\eq c\1\ee_1$ to the above equation (which becomes $\nu_1\pl 2\1\msum_{i=2}^n\,\nu_i\eq ab$).
\sk
In order to understand the above calculation better, it may be good to consider the projection $\pi_1:\R^n\tos\R^{n-1}$ (which forgets the first coordinate). Indeed, we must have the inequality
$$|\pi_1(a\1\ob)|\eq a(n{-}1)\sgt|\pi_1(b\1\ff_i)|\eq b(n{-}2),$$
\par\nin since the first coordinate of $a\1\ob$ is smaller than that of $b\1\ff_i$. Here one can also consider the projection $\R^n\tos\R^2$ defined by $(\nu_i)\mapsto\bl(\nu_1,\sum_{i=2}^n\nu_i\br)$. It seems easier to calculate the equation of the line containing the image of the face.
\sk
Assume now $n\eq 4$ and moreover the vertices of $\Gp(f)$ consist of $\,a\ob$, $b\1\ff_i$, $c\1\ee_i$ ($i\ins[1,n])$ for simplicity. Then we can verify that $(a,b,c)\eq(5,7,36)$, where the numbers of $k$-dimensional compact faces of $\Gp(f)$ are $9,28,34,14$ for $k\eq0,1,2,3$ respectively (with $9\mi 28\pl 34\mi 14\eq 1$ and $\mu_f\eq 385437$). In this case it does not seem quite easy to draw a picture.\par
$$\setlength{\unitlength}{.5cm}
\begin{picture}(12,10.5)
\multiput(4,4)(.2,0){35}{\circle*{.1}}
\multiput(4,4)(0,0.2){30}{\circle*{.1}}
\multiput(1,0)(.12,.16){25}{\circle*{.1}}
\multiput(1,0)(.06,.2){50}{\circle*{.1}}
\multiput(1,0)(.2,.08){50}{\circle*{.1}}
\multiput(4,10)(.14,-.12){50}{\circle*{.1}}
\qbezier(3.1,4.3)(2.05,2.15)(1,0)
\qbezier(3.1,4.3)(3.55,4.15)(4,4)
\qbezier(3.1,4.3)(3.55,7.15)(4,10)
\qbezier(4.8,3)(2.9,1.5)(1,0)
\qbezier(4.8,3)(4.4,3.5)(4,4)
\qbezier(4.8,3)(7.9,3.5)(11,4)
\qbezier(5.9,5)(3.45,2.5)(1,0)
\qbezier(5.9,5)(4.95,7.5)(4,10)
\qbezier(5.9,5)(8.35,4.5)(11,4)
\qbezier(6,6)(5,5)(4,4)
\qbezier(6,6)(5,8)(4,10)
\qbezier(6,6)(8.5,5)(11,4)
\linethickness{0.4mm}
\qbezier(3.1,4.3)(3.95,3.65)(4.8,3)
\qbezier(3.1,4.3)(4.5,4.65)(5.9,5)
\qbezier(3.1,4.3)(4.55,5.15)(6,6)
\qbezier(4.8,3)(5.35,4)(5.9,5)
\qbezier(4.8,3)(5.4,4.5)(6,6)
\qbezier(5.9,5)(5.95,5.5)(6,6)
\qbezier(4.8,4.3)(3.95,4.3)(3.1,4.3)
\qbezier(4.8,4.3)(4.8,3.65)(4.8,3)
\qbezier(4.8,4.3)(5.35,4.65)(5.9,5)
\qbezier(4.8,4.3)(5.4,5.15)(6,6)
\put(1,0){\circle*{.23}}
\put(.25,0){$\scriptstyle r_1$}
\put(4,4){\circle*{.23}}
\put(3.03,3.67){$\scriptstyle r_4$}
\put(4,10){\circle*{.23}}
\put(3.2,9.9){$\scriptstyle r_3$}
\put(11,4){\circle*{.23}}
\put(11.3,3.9){$\scriptstyle r_2$}
\put(4.8,4.3){\circle*{.23}}
\put(5.4,3.6){$\scriptstyle p_0$}
\put(3.1,4.3){\circle*{.23}}
\put(2.38,4.3){$\scriptstyle p_2$}
\put(4.8,3){\circle*{.23}}
\put(4.8,2.5){$\scriptstyle p_3$}
\put(5.9,5){\circle*{.23}}
\put(5.9,4.5){$\scriptstyle p_4$}
\put(6,6){\circle*{.23}}
\put(6.2,6.2){$\scriptstyle p_1$}
\end{picture}$$
\par\nin \skn
This is a generalization of the picture in Example~\hl{E2.2}{2.2} just above, which is obtained by taking the intersections of cones of 0 or 1-dimensional compact faces with the hyperplane defined by $\msum_{i=1}^4\,\nu_i\eq1$, where the coordinate axes as in the pictures of Example~\hl{E2.2}{2.2} are omitted. The points corresponding to $a\1\ob$, $b\1\ff_i$, and $c\1\ee_i$ are denoted respectively by $p_0$, $p_i$, and $r_i$ ($i\eq 1,\dots,4$). The segments corresponding to 1-dimensional compact faces $\si$ with $k(\si)\eq 4$, 3, and 2 are denoted respectively by thick, thin, and dotted lines. These consist respectively of
$$\begin{array}{lll}{[p_i,p_j]}&\bl(i,j\ins[0,4],\,i\slt j\br)&{[4;\,10],}\\{[p_i,r_j]}&\raise5mm\h{}\bl(i,j\ins[1,4],\,i\ne j\br)&{[3;\,12],}\\{[r_i,r_j]}&\raise5mm\h{}\bl(i,j\ins[1,4],\,i\slt j\br)&{[2;\,6].}\end{array}$$ 
where $[p_i,p_j]$ denotes the convex hull of $\{p_i,p_j\}$, etc., and similarly for $[p_i,p_j,p_k]$, etc. The 2 and 3-dimensional compact faces of $\Gp(f)$ correspond respectively to
$$\begin{array}{lll}{[p_i,p_j,p_k]}&\bl(i,j,k\ins[0,4],\,i\slt j\slt k\br)&{[4;\,10],}\\{[p_i,p_j,r_k]}&\raise5mm\h{}\bl(i,j,k\ins[1,4],\,k\,{\ne}\,i\slt j\,{\ne}\,k\br)&{[4;\,12],}\\{[r_i,r_j,p_k]}&\raise5mm\h{}\bl(i,j,k\ins[1,4],\,k\,{\ne}\,i\slt j\,{\ne}\,k\br)&{[3;\,12],}\end{array}$$ 
and
$$\begin{array}{lll}{[p_0,p_i,p_j,p_k]}&\bl(i,j,k\ins[1,4],\,i\slt j\slt k\br)&{[4;\,4],}\\{[p_i,p_j,p_k,r_l]}&\raise5mm\h{}\bl(\{i,j,k,l\}\eq\{1,\dots,4\},\,i\slt j\slt k\br)&{[4;\,4],}\\{[p_i,p_j,r_k,r_l]}&\raise5mm\h{}\bl(\{i,j,k,l\}\eq\{1,\dots,4\},\,i\slt j,\,k\slt l\br)&{[4;\,6].}\end{array}$$ 
At the end of each line, the value of $k(\si)$ and the number of corresponding faces are written. (Note that $10\pl 12\pl 6\eq 28$, etc.)
We can then calculate $\ep_{\si},\ep_{\tau}$ in Theorem\,\,\hl{T4}{4} as follows:
$$\aligned\ep_{\si}&\eq 3\q\h{if}\,\,\,\si\eq[p_i]\,\,\,(i\ins[1,4]),\\
\ep_{\tau}&\,{=}\begin{cases}2&\h{if}\,\,\,\tau\eq[p_i,p_j]\,\,\,(i,j\ins[1,4],\,i\slt j),\\ 1&\h{if}\,\,\,\tau\eq[p_i,r_j]\,\,\,(i,j\ins[1,4],\,i\nes j).\end{cases}\endaligned$$
\par\nin These vanish otherwise.
\par\htt{R2.4a}{}\msn
{\bf Remark\,\,2.4a.} When $n\eq 3$, we can verify that $(a,b,c)\eq(3,5,16)$ with $\mu_f\eq 992$.
\par\htt{R2.4b}{}\msn
{\bf Remark\,\,2.4b.} We may consider a similar quiz replacing the $\ff_i$ with $\ee_{i,j}\defs\ee_i\pl\ee_j$ ($i\slt j$), where the condition on the maximal-dimensional compact face is about $a\1\ob$, $b\1\ee_{1,j}$ ($j\ins[2,n]$). Here we assume $n\gess4$, since it is the same as above when $n\eq 3$. The case $n\gess5$ may be rather complicated. Assume $n\eq 4$.
The minimal triple then seems to be given by $(a,b,c)\eq(2,5,21)$. We have the following picture.\par
$$\setlength{\unitlength}{.5cm}
\begin{picture}(12,10.5)
\multiput(4,4)(.2,0){35}{\circle*{.1}}
\multiput(4,4)(0,0.2){30}{\circle*{.1}}
\multiput(1,0)(.12,.16){25}{\circle*{.1}}
\multiput(1,0)(.06,.2){50}{\circle*{.1}}
\multiput(1,0)(.2,.08){50}{\circle*{.1}}
\multiput(4,10)(.14,-.12){50}{\circle*{.1}}
\put(1,0){\circle*{.23}}
\put(.25,0){$\scriptstyle r_1$}
\put(4,4){\circle*{.23}}
\put(3.37,4.28){$\scriptstyle r_4$}
\put(4,10){\circle*{.23}}
\put(3.2,9.9){$\scriptstyle r_3$}
\put(11,4){\circle*{.23}}
\put(11.3,3.9){$\scriptstyle r_2$}
\put(4.8,4.3){\circle*{.23}}
\put(2.5,2){\circle*{.23}}
\put(2.5,5){\circle*{.23}}
\put(6,2){\circle*{.23}}
\put(4,7){\circle*{.23}}
\put(7.5,4){\circle*{.23}}
\put(7.5,7){\circle*{.23}}
\qbezier(2.5,2)(2.5,3.5)(2.5,5)
\qbezier(2.5,2)(4.25,2)(6,2)
\qbezier(2.5,2)(3.25,4.5)(4,7)
\qbezier(2.5,2)(5,3)(7.5,4)
\qbezier(2.5,5)(4.25,3.5)(6,2)
\qbezier(2.5,5)(3.25,6)(4,7)
\qbezier(2.5,5)(5,6)(7.5,7)
\qbezier(6,2)(6.75,3)(7.5,4)
\qbezier(6,2)(6.75,4.5)(7.5,7)
\qbezier(4,7)(5.75,5.5)(7.5,4)
\qbezier(4,7)(5.75,7)(7.5,7)
\qbezier(7.5,4)(7.5,5.5)(7.5,7)
\linethickness{0.4mm}
\qbezier(4.8,4.3)(3.65,3.15)(2.5,2)
\qbezier(4.8,4.3)(3.65,4.65)(2.5,5)
\qbezier(4.8,4.3)(5.4,3.15)(6,2)
\qbezier(4.8,4.3)(4.4,5.65)(4,7)
\qbezier(4.8,4.3)(6.15,4.15)(7.5,4)
\qbezier(4.8,4.3)(6.15,5.65)(7.5,7)
\end{picture}$$
\par\nin \bs\bs\htt{S3}{}
\vbox{\centerline{\bf 3. Simplicial Newton polytope case}
\bsn
In this section we prove the main theorems after introducing $\Ga$-spectrum and explaining the Steenbrink conjecture on spectral pairs.}
\par\htt{3.1}{}\msn
{\bf 3.1.~$\Ga$-spectrum.} In this paper we say that a polynomial $f$ of $n$ variables has \h{\it simplicial\1} {\it Newton boundary\1} or the Newton polytope $\Gp(f)$ is {\it simplicial\1} (or $f$ is simplicial for short) if any {\it compact\1} face $\si$ of $\Gp(f)$ is a simplex, that is, $\si$ is a convex hull of $k{+}1$ vertices with $k\defs d_{\si}\,({=}\,\dim\si)$. Note that $\Gp(f)$ is always simplicial when $n\eq2$.
\sk
Assume $f$ is non-degenerate and simplicial. Let $\si$ be an $r$-dimensional face of $\Gp(f)$ which is the convex hull of the vertices $\xi^{(k)}\ins\N^n$ ($k\ins[0,r]$). Set
$$\aligned&g_{\si}\defs\mprod_{k=0}^n\,(1\pl x^{\xi^{(k)}}),\q E_{\si}\defs(\Spx g_{\si})^{\rm con.hull},\\&\q\q\q\q q_{\si}(t)\defs\msum_{\nu\in\N^n\cap\1E_{\si}^{\circ}}\,t\1^{\ell(\nu)}.\endaligned$$
\par\nin Here $S\1^{\rm con.hull}$ denotes the convex hull of a subset $S\sst\R^n$ in general, $E_{\si}^{\circ}$ is the interior of $E_{\si}$ in the affine space spanned by $E_{\si}$, and $\ell$ is any linear function without a constant term such that $\ell^{-1}(1)\supset\si$.
\sk
We say that a face $\si\sst\Gp(f)$ is {\it interior,} if it is not contained in any coordinate hyperplane. Let $\CFfin$ be the set of interior compact faces of $\Gp(f)$, and $\CFfin^k$ the subset of $k$-dimensional faces. Define the $\Ga$-{\it spectrum\1} of $f$ by
\htt{3.1.1}{}
$$\aligned\SpG_f(t)&=\msum_{\si\in\CFfin}\,\bl(\msum_{j=0}^{c_{\si}}\,t^j\br)q_{\si}(t)\\&\q +\bl|\CFfin^0\br|\,\bl(\msum_{j=0}^{n-2}\,t^j\br)\1t.\endaligned
\leqno(3.1.1)$$
\par\nin Here $c_{\si}\defs n{-}1{-}d_{\si}$, and $|S|$ denotes the number of elements of a finite set $S$ in general.
\sk
In the isolated singularity case, the multiplication by $t$ corresponds to the action of $\Gr_Vf$ on $\Gr_V^{\ssb}\bl(\C\{x\}/(\dd f)\br)$, and to that of $N$ on $\Gr_F^{\ssb}H^{n-1}(\Ff,\C)$. Hence the summation {\it without the multiplicative factors\1} $\bl(\msum_{j=0}^{c_{\si}}\,t^j\br)$, $\bl(\msum_{j=0}^{n-2}\,t^j\br)$ is called the $N$-{\it primitive part\1} of the $\Ga$-spectrum. The latter corresponds to part of the $N$-primitive part of $\Gr^W_{\ssb}H^{n-1}(\Ff,\C)$ by the modified Steenbrink conjecture for spectral pairs.
\sk
Setting
$$\SpGp_f(t)\defs\SpG_f(t^{-1})\,t^n,$$
\par\nin we have the {\it symmetry}
\htt{3.1.2}{}
$$\SpG_f(t)\eq\SpGp_f(t),
\leqno(3.1.2)$$
\par\nin since $q_{\si}(t)\eq q_{\si}(t^{-1})\,t^{d_{\si}+1}$ by the definition of $q_{\si}(t)$. Set
$$\SpDef_f(t)\defs\Sp_f(t)-\SpG_f(t),\q\Sp_f^{\prime\,{\rm Def}}(t)\defs\Sp'_f(t)-\SpG_f(t).$$
\par\nin These are called the {\it defects.} In the isolated singularity case, we have the symmetry
\htt{3.1.3}{}
$$\SpDef_f(t)\eq\Sp_f^{\prime\,{\rm Def}}(t),
\leqno(3.1.3)$$
\par\nin since $\Sp_f(t)\eq\Sp'_f(t)$. The spectral numbers belonging to these are called {\it defective\1} ones. By the definition of Newton filtration $V_N^{\ssb}$ (see (\hl{1.2.4}{1.2.4})), we have
\htt{3.1.4}{}
$$(\dd f)\sst V_N^{>1}\C\{x\}.
\leqno(3.1.4)$$
\par\nin Comparing this with the definition of $\Ga$-spectrum, we see that the defective spectral numbers are greater than 1, since the Newton filtration $V_N^{\al}\C\{x\}$ is generated by monomials ($\al\ins\Q$), see also a remark after the main theorem of \cite{exp}. In view of the self-symmetry (\hl{3.1.3}{3.1.3}), we then get the following.
\par\htt{L3.1}{}\msn\vbox{\nin
{\bf Lemma\,\,3.1.} {\it If $f$ is non-degenerate and convenient, then
\htt{3.1.5}{}
$$\{\h{Defective spectral numbers}\}\sst(1,n{-}1).
\leqno(3.1.5)$$
\par\nin In particular, for $n\eq2$, we have $\SpDef_f(t)\eq0$, that is,}
\htt{3.1.6}{}
$$\Sp_f(t)=\msum_{\si\in\CFfin}\,\bl(\msum_{j=0}^{c_{\si}}\,t^j\br)q_{\si}(t)\,+\bl|\CFfin^0\br|\,t.
\leqno(3.1.6)$$}
\sk
(Note that (\hl{3.1.6}{3.1.6}) also follows from \cite{St1}, see \cite{Ar} and the picture in Example~\hl{E1.4b}{1.4b}.)
\sk
As a corollary of Proposition\,\,\hl{P1.4}{1.4}, we then get the following
\par\htt{P3.1}{}\msn
{\bf Proposition\,\,3.1.} {\it Assume $n\eq2$, and $f$ is non-degenerate and non-convenient. For $i\ins I_f$, let $\si_i$ be the $0$-dimensional face of $\Gp(f)$ contained in the non-compact face parallel to the $x_i$-axis, where $|I_f|\eq1$ or $2$ by the non-convenience condition. Then}
\htt{3.1.7}{}
$$\Sp'_f(t)\eq\SpG_f(t)-\msum_{i\ins I_f}\,q_{\si_i}(t).
\leqno(3.1.7)$$
\par\nin \ms
(Note that $q_{\si_i}(t)\eq0$ for $i\ins I_f$ if $f$ has an isolated singularity at 0.)
\ms
When $n\eq 3$, the defective spectral numbers are contained in the open interval $(1,2)$, and can be determined by comparing the definition of $\Ga$-spectrum with the formula for monodromy zeta functions \cite{Va1}. This can give a proof of Theorem\,\,\hl{T1}{1}.
\par\htt{R3.1a}{}\msn
{\bf Remark\,\,3.1a.} For a polynomial $f$ of 3 variables with a simplicial convenient non-degenerate Newton boundary, one can calculate $\SpG_f(t)$ from the Newton polyhedron $\Gp(f)$ using a computer. Combining this with a formula for monodromy zeta functions \cite{Va1}, one can determine the spectrum $\Sp_f(t)$, since the defective spectral numbers are contained in $(1,2)$. It is also possible to compute $\Sp_f(t)$ using (\hl{2.1.3}{2.1.3}). These two calculations agree as far as computed by writing a small computer program, and are compatible with a computation using Singular \cite{Sing} as far as examined. (Note that the spectral numbers are shifted by $-1$ in the latter as in \cite{St2}.)
\par\htt{R3.1b}{}\msn
{\bf Remark\,\,3.1b.} For $n\eq 3$, Theorem\,\,\hl{T1}{1} implies that the defect $\SpDef_f(t)$ is a linear combination of fractional power polynomials
$$t\1\bl(\msum_{i=1}^{e_j-1}\,t\1^{i/e_j}\br),$$
\par\nin with $e_j$ the greatest common divisor of the coordinates of {\it certain\1} vertices $P_j\ins\Gp(f)$. It seems that the monomials corresponding to defective spectral numbers are given by the points of
$$\mcup_j(Q_j,P_j{+}Q_j)\cap\N^n.$$
\par\nin Here the $Q_j$ are {\it certain\1} vertices of $\Gp(f)$ such that the convex hull of $\{P_j,Q_j\}$ is a face of $\Gp(f)$, and $(Q_j,P_j{+}Q_j)$ denotes the {\it interior\1} of the convex hull of $\{Q_j,P_j{+}Q_j\}$ in the line passing through $Q_j$, $P_j{+}Q_j$. It seems, however, very difficult to determine which vertices are the $Q_j$, since their information is not sufficiently reflected to numerical data such as defective spectral numbers. (They are partially encoded in the combinatorial polynomials $r_{\!\si}(t)$ in (\hl{2.1.3}{2.1.3}), although some cancelations may occur.)
\par\htt{R3.1c}{}\msn
{\bf Remark\,\,3.1c.} In the {\it non-simplicial\1} case, Theorem\,\,\hl{T1}{1} and the properties (\hl{1}{1}) and (\hl{2}{2}) do not hold even though the $q_{\si}(t)$ can be defined as in \cite{St1}. For instance, setting
$$f\eq x^2\pl y^2\pl xz\pl yz\pl z^4,$$\vskip-3mm\nin
we have
$$\SpG_f(t)\eq q_{\si_1}\eq t\pl t^{3/2},\q\h{although}\q\Sp_f(t)\eq t^{3/2},$$
\par\nin see \cite{Da}. Here $\si_1$ is the unique {\it non-simplicial\1} 2-dimensional face of $\Gp(f)$. More precisely, the other $q_{\si}(t)$ vanish except for $q_{\emptyset}(t)\eq1$, and the combinatorial polynomials for these are given by $r_{\si_1}(t)\eq1$, $r_{\emptyset}(t)\eq{-}t$. This shows that Steenbrink conjecture on spectral pairs does not hold in the {\it non-simplicial\1} case.
(This seems to be related to the non-triviality of the {\it primitive\1} cohomology of the toric variety associated with a {\it non-simplicial\1} face.)
\par\htt{3.2}{}\msn
{\bf 3.2.~Spectral pairs.} The example in Remark\,\,\hl{R3.1c}{3.1c} was found as a counter-example to Steenbrink's conjecture on spectral pairs \cite{St1} in the {\it non-simplicial\1} case, see \cite{Da}. Recall that the {\it spectral pairs\1} consist of $\mu_f$ pairs of rational numbers $(\al_i,w_i)$ ($i\ins[1,\mu_f]$) with $\al_i$ the spectral numbers and $w_i$ the modified weights. In this paper, the weights are shifted by $-1$ when $\al_i\ins\Z$. So the center of symmetry of the weights is always $n{-}1$. More precisely, we define the spectral pairs in this paper as follows (using the symmetry of spectral numbers, see \hl{1.1}{1.1} and \cite[Sections 2.1--2]{ste}):
$$\aligned&\#\{\1i\mid(\al_i,w_i)\eq(\al,w)\}=\dim\Gr_F^p\Gr^W_{w+\de_{\la,1}}H^{n-1}(\Ff,\C)_{\la}\\ &\q\q\h{with}\q\q p\eq[\al],\q\la\eq\exp(2\pi\sqrt{-1}\1\al),\endaligned$$
\par\nin where $\de_{\la,1}\eq1$ if $\la\eq1$, and 0 otherwise.
\sk
This definition essentially coincides with the one in \cite{St1}, and is different from the one in \cite{SSS} by the shift of spectral numbers by $1$ and also by the change of modified weights under the involution of $\Z$ defined by $w\mapsto 2n\mi2\mi w$, where the {\it symmetry of modified weights\1} with center $n{-}1$ is used. Here it is actually better to change the spectral numbers (instead of modified weights) applying the involution of $\Q$ defined by $\al\mapsto n\mi\al$. Here the {\it self-duality\1} of the mixed Hodge structure on the vanishing cohomology is used.
\sk
In the notation of \hl{3.1}{3.1}, we define the {\it weighted spectrum\1} $\SpW_f(t,u)$ to be the generating polynomial of the spectral pairs, that is,
$$\SpW_f(t,u):=\msum_{i=1}^{\mu_f}\,t^{\al_i}u^{w_i}.$$
\par\nin We can define also the {\it weighted $\Ga$-spectrum\1} by
    $$\aligned\SpWG_f(t,u):={}&\msum_{\si\in\CFfin}\,\bl(\msum_{j=0}^{c_{\si}}\,t^ju^{2j}\br)q_{\si}(t)u^{d_{\si}}\\{}+{}&\bl|\CFfin^0\br|\,\bl(\msum_{j=0}^{n-2}\,t^ju^{2j}\br)\1tu,\endaligned$$
\par\nin together with the weighted defect
$$\SpDefW_f(t,u):=\SpW_f(t,u)\mi\SpWG_f(t,u).$$
\par\nin \par\htt{R3.2a}{}\msn
{\bf Remark\,\,3.2a.} Steenbrink's conjecture predicted that the spectral pairs could be obtained by adding the information of the modified weights to the last term of (\hl{5}{5}), see \cite{St1}. Here the problem is the modified weight of the trivial Poincar\'e series $q_{\emptyset}(t)\eq1$ for $\si\eq\emptyset$\,: it must be $-1$, since the modified weight of $t\1q_{\emptyset}(t)$ and $t^{n-1}q_{\emptyset}(t)$ should be 1 and $2n{-}3$ respectively; for instance, if
$$f\eq x^4\pl y^4\pl z^4\pl xyz.$$
\par\nin Note that the multiplication by $t$ corresponds to that of $\Gr_Vf$, and further to the action of $-N/2\pi i$ on the vanishing cohomology, see a remark after (\hl{1.1.7}{1.1.7}). The latter is a morphism of mixed Hodge structures of type $(-1,-1)$ so that it {\it decreases\1} the weights by $2$.
This also holds for the normalization of weights in \cite{SSS} explained above. With our normalization, however, the action of $t$ {\it increases\1} the weights by 2 because of the {\it involution\1} of weights explained above.
This comes from the definition of spectral pairs written above, since we adopt the definition of dual (or Hodge) spectrum $\Sp'_f(t)$ as in (\hl{1.1.3}{1.1.3}).
\sk
In the case of the example in Remark\,\,\hl{R3.1c}{3.1c}, the modified weight of $q_{\emptyset}(t)\eq1$ must be 0 in order that $-t\1q_{\emptyset}(t)\eq{-}t$ has the {\it middle\1} weight 2, and cancels out with $t$ in $\SpG_f(t)$.
\par\htt{R3.2b}{}\msn
{\bf Remark\,\,3.2b.} In the {\it simplicial\1} case, the above problem does not occur by (\hl{5}{5}) together with the properties (\hl{1}{1}) and (\hl{2}{2}) in the introduction. The weight of $t^jq_{\si}(t)$ should be given by
\htt{3.2.1}{}
$$d_{\si}\pl 2j\q\q\bl(\,j\ins[0,n{-}1{-}d_{\si}]\,\br),
\leqno(3.2.1)$$
\par\nin independently of $k(\si)$ so that a modified {\it Steenbrink conjecture\1} in the simplicial convenient case can be formulated as
\htt{3.2.2}{}
$$\SpW_f(t,u)=\msum_{\emptyset\les\si\les\Gf}\,r_{\!\si}(tu^2)\1q_{\si}(t)\1u^{d_{\si}},
\leqno(3.2.2)$$
\par\nin with $r_{\!\si}(t)$ as in (\hl{5}{5}) (and $d_{\si}\eq d(\si)\mi1$). This is a consequence of the {\it descent theorem of nearby cycle formula,} see \cite[Theorem\,\,2]{des}. We have the {\it symmetry\1} of combinatorial polynomials as in (\hl{6}{6}) (that is, Theorem\,\,\hl{TA}{A} in Appendix). This is compatible with the {\it monodromical property\1} of the weight filtration in (\hl{1.1.8}{1.1.8}).
\sk
By the proof of Theorem\,\,\hl{T1}{1} in \hl{3.3}{3.3} below, the coefficient of $q_{\si}$ in (\hl{4}{4}) coincides with $r_{\!\si}(t)$. Here $q_{\emptyset}(t)\eq1$ so that $r_{\emptyset}(t)$ coincides with the second term of the right-hand side of (\hl{4}{4}). This implies Theorem\,\,\hl{T2}{2}.
\par\htt{R3.2c}{}\msn
{\bf Remark\,\,3.2c.} Let $n_{\la,k}$ be the number of Jordan blocks of the Milnor monodromy for the eigenvalue $\la$ with size $k$. The formula (\hl{3.2.2}{3.2.2}) is closely related to these numbers. Let $\de_{\si}$ be the positive integer such that
$$\ell_{\si}(\Z^n\cap V_{\si})\eq\tfrac{1\,\,}{\de_{\si}}\,\Z,$$
\par\nin where $V_{\si}\sst\R^n$ is the $\R$-vector subspace spanned by $\si$, and $\ell_{\si}$ is the linear function on $V_{\si}$ with $\si\sst\ell_{\si}^{-1}(1)$. (In some literature, $\de_{\si}$ may be called the {\it lattice distance\1} from the origin.) Set
$$\aligned\CFfla^k&\defs \{\si\ins\CFf^k\mid\la^{\de_{\si}}\eq1\},\\ \CFfinla^k&\defs\CFfin^k\cap\CFfla^k.\endaligned$$
\par\nin \sk
It seems that the following is known to specialists, see for instance \cite{Sta} and references there (although some arguments using Hodge theory do not seem quite easy to follow, since the proof of the descent theorem of nearby cycle formula is not entirely trivial, see also \cite{des}):
\htt{3.2.3}{}
$$n_{\la,k}=\begin{cases}\,\bl|\CFfinla^0\br|&(\la{\ne}1,\,k\eq n),\\ \,\msum_{\si\in\CFfinla^1}\,l(\si)\,-\,\msum_{\si\in\CFfla^0}\,\be^1_{\si}\raise5mm\h{}&(\la{\ne}1,\,k\eq n{-}1),\\ \,\bl|\,\mcup_{\si\in\CFfin^1}\1\si\cap\Z_{>0}^n\1\br|\raise5mm\h{}&(\la\eq1,\,k\eq n{-}1).\end{cases}
\leqno(3.2.3)$$
\par\nin Here $l(\si)\defs|\Z^n\cap\si|-1$ for $\si\ins\CFf^1$, and $\be^1_{\si}$ is the number of 1-dimensional interior compact faces of $\Gp(f)$ containing $\si\ins\CFf^0$. It is well known that $n_{1,n}\eq0$. We can deduce (\hl{3.2.3}{3.2.3}) from (\hl{3.2.2}{3.2.2}). Note that (\hl{3.2.3}{3.2.3}) (especially for $k\eq n$) is a direct consequence of the descent theorem of nearby cycle formula using a standard estimate of weights \cite[4.5.2]{mhm}, where the theory of Danilov and {Khovanski\u\i} is not necessary, see also \cite{des}.
\par\htt{R3.2d}{}\msn
{\bf Remark\,\,3.2d.} The {\it Thom-Sebastiani theorem for weighted spectrum\1} is expressed by
\htt{3.2.4}{}
$$\SpW_{f+g}(t,u)\eq\SpW_f(t,u)\,\SpW_g(t,u)\1u,
\leqno(3.2.4)$$
\par\nin with multiplication taken in $\Z[t^{1/e},u]$ for a sufficiently divisible positive integer $e$. Here $f,g$ are holomorphic function on complex manifolds $X,Y$ with isolated singularities, and $f\pl g$ is a function on $X{\times}\1Y$.
We need the multiplication by $u$ on the right-hand side of (\hl{3.2.4}{3.2.4}), since the middle weight is given by
$$d_{X\times Y}{-}1\eq(d_X{-}1)\pl(d_{Y}{-}1)\pl 1.$$
\par\nin \sk
There seems to be an error in the formulation of the theorem in \cite{SS} as is noted in \cite{KL1}. However, this is never serious, since it occurred during the conversion between two systems of weights of spectral pairs: one is as in this paper, and the other shifts the weight by 1 on the unipotent monodromy part as in \cite{KL1}. This conversion is given by adding $\langle\al_i\rangle$ to the weight $w_i$ for a spectral pair $(\al_i,w_i)$. Here $\langle\al\rangle\eq 1$ if $\al\ins\Z$, and 0 otherwise. This produces the {\it correction term\1} for weights in the Thom-Sebastiani isomorphism defined by
$$\langle\al\,|\,\be\rangle\defs 1\,{+}\1\langle\al\pl\be\rangle\mi\langle\al\rangle\mi\langle\be\rangle,$$
\par\nin see also \cite{KL1}. Note that the monodromy filtration is {\it compatible\1} with the tensor product of two vector spaces endowed with a nilpotent endomorphism $N$.
(Recall that the Hodge and weight filtrations $F, W$ were treated separately in \cite{SS}. For $W$, its coincidence with the monodromy filtration up to shift was used. For $F$, the Brieskorn lattice together with Malgrange's idea was hired.) The above formula (\hl{3.2.4}{3.2.4}) is compatible with the involution of weights describing the difference from the definition in \cite{SSS}, which is explained just after the definition of spectral numbers. (This becomes clearer if we shift the weights so that the center of symmetry is always 0.)
\par\htt{3.3}{}\msn
{\bf 3.3.~Proofs of Theorems~\hl{T1}{1}--\hl{T2}{2} and \hl{T4}{4}--\hl{T5}{5}.} The assertions for spectrum follow from the Steenbrink formula for spectrum (\hl{5}{5}) (with \cite{exp}) by using the symmetry of combinatorial polynomials (\hl{6}{6}) and also the estimate of degree (\hl{7}{7}). Indeed, these imply that $r^{\rm Def}_{\!\si}(t)\eq 0$ if $d_{\si}\ne0$ ($n\eq3$) or $d_{\si}\ne-1,0,1$ ($n\eq4$).
We can get (\hl{3}{3}) and (\hl{11}{11}) from the definition of $r_{\!\si}(t)$ written in (\hl{5}{5}) by calculating only the coefficients of higher degree terms.
By the definition of $r^{\rm Def}_{\!\si}(t)$ written after (\hl{6}{6}), we get $\ep_{\si}$ (or $\ep_{\tau}$) and $\ep_{\emptyset}$ in (\hl{11}{11}) by subtracting $1$ and $n^{0,4}$ from the corresponding coefficient of $r_{\!\si}(t)$ (or $r_{\tau}(t)$) and $r_{\emptyset}(t)$ respectively if $\si,\tau$ are interior, that is, if $k(\si)\eq k(\tau)\eq 4$. This is due to the difference from $\SpG_f(t)$. The argument is similar for (\hl{3}{3}). We can get a formula similar to the second equality of (\hl{11}{11}), see also Remark\,\,\hl{R3.3}{3.3} below.
\sk
The assertions for weighted spectrum then follow from the modified Steenbrink conjecture for spectral pairs (see \cite[Theorem\,\,2]{des} and (\hl{3.2.2}{3.2.2})) as is noted in Remark\,\,\hl{R3.2b}{3.2b} for $n\eq3$. This finishes the proofs of Theorems~\hl{T1}{1}, \hl{T2}{2} and \hl{T4}{4}.
\par\htt{R3.3}{}\msn
{\bf Remark\,\,3.3.} The assertions related to $\ga_{\si}$, $\ga_{\tau}$ in Theorems~\hl{T1}{1} and \hl{T4}{4} can be generalized to the one about the coefficient $r^{\rm Def}_{\si,1}$ of $\,t\,$ in $\,r^{\rm Def}_{\!\si}(t)\eq r^{\rm Def}_{\si,1}t\,$ for $\si\ins\CFf^{n-3}$.
This is in terms of the number of (not necessarily compact) $(n{-}1)$-dimensional faces of $\Gp(f)$ containing $\si\ins\CFf^{n-3}$, and can be shown by restricting to a transversal subspace intersecting $\si$ at one point.
\par\htt{3.4}{}\msn
{\bf 3.4.~Proof of Theorem\,\,\hl{T3}{3}.} For $r\gg 0$, set
$$\CFflc^j\defs\CFfl^j\stm\CFf,$$
\par\nin and similarly for $\CFflinc^j$. (In general $\,_{\bf c}\,$ denotes the {\it complement\1} of $\CFf$.) Note that $\CFflinc^0\eq\emptyset$. By \cite[Lemma\,\,A.2]{JKSY}, there is a canonical decomposition indexed by $I_f$
$$\CFflc^j=\h{$\bigsqcup$}_{i\in I_f}\CFflci^j\,,$$
\par\nin together with the canonical bijections
\htt{3.4.1}{}
$$\CFfbi^j\simto\CFflinci^{j+1}\q(j\eq0,1),
\leqno(3.4.1)$$
\par\nin sending $\tau$ to the convex hull of $\si^{(i)}\cup\tau$ with $\si^{(i)}\defs\{r\1\ee_i\}$ (and $\CFflci^0\eq\{\si^{(i)}\}$). Here $\CFfbi^j\sst\CFf^j$ ($i\ins I_f$) is the subset consisting of $\tau$ such that the convex hull of $\si^{(i)}\cup\tau$ is a face of $\Gp(f{+}\ell\1^r)$ ($r\gg 0$). Set
$$\CFfb^j=\mcup_{i\in I_f}\CFfbi^j,\q\CFfb=\h{$\bigsqcup$}_j\,\CFfb^j.$$
\par\nin We have $m_{\tau}\eq0$ if $\tau\notin\CFfb^j$. However, the converse does not necessarily hold because of the last condition in the definition of $m_{\si}$, that is, $\si^{(i)}\cup\tau$ is not contained in any coordinate plane (consider the case where $j\eq0$ and $\tau$ is contained in a coordinate plane containing $\si^{(i)}$). Note also that the above condition on $\si^{(i)}\cup\tau$ can be different from that $\tau$ is not contained in any coordinate plane, since $\tau$ may be contained in the coordinate plane {\it not containing\1} $\si^{(i)}$ (for instance, if $f\ins\C[x,y]\sst\C[x,y,z]$).
\sk
We compare the part {\it depending really on} $r$ of the right-hand side of (\hl{4}{4}) in Theorem\,\,\hl{T1}{1} applied to $f{+}\ell\1^r$ with that of (\hl{1.3.1}{1.3.1}). In view of Proposition\,\,\hl{P1.4}{1.4} and Lemma\,\,\hl{L1.4}{1.4}, we then see that the right-hand side of (\hl{1.3.1}{1.3.1}) can be expressed by
\htt{3.4.2}{}
$$\aligned&\msum_{\si\in\CFflinc}\,\bl(\msum_{j=0}^{2-d_{\si}}\,t^j\br)q_{\si}(t)\\
&+\msum_{\si\in\CFflc^0}(\gat_{\si}{-}3)\1t\1q_{\si}(t)\raise12pt\h{}\\
&+\msum_{\tau\in\CFfb}\,m_{\tau}\bl(\msum_{j=0}^{1-d_{\tau}}\,t^j\br)q_{\tau}(t)+\msum_{\tau\in\CFfb^0}\,m_{\tau}\1t.
\endaligned
\leqno(3.4.2)$$
\par\nin using the restriction to the summation over the $(i,j)$ with $\al_{i,j}\less 1$, see Remark\,\,\hl{R1.3}{1.3}. Here we first determine which part of (\hl{4}{4}) {\it really\1} depends on $r$, and then find an {\it appropriate term\1} which is independent of $r$ and whose sum with the above term depending on $r$ coincides with the right-hand side of (\hl{3.4.2}{3.4.2}). We first verify that the first two terms of (\hl{3.4.2}{3.4.2}) {\it really\1} depend on $r$ (since $q_{\si}(t)$ is defined by using the lattice points in the {\it interior\1} of the parallelotope $E_{\si}$, and $\si$ is the convex hull of $\{r\1\ee_i\}\cup\tau$ for $\tau\ins\CFfbi^j$). The last two terms of (\hl{3.4.2}{3.4.2}) give a generalization of the circled point in the picture of Example~\hl{E1.4b}{1.4b} (with coordinates $(1,0)$). They can be written as the sum over $i\ins I_f$ of
\htt{3.4.3}{}
$$\msum_{\tau\in\CFfbi}\,\bl(\msum_{j=0}^{1-d_{\tau}}\,t^j\br)q_{\tau}(t)+\bl(\bl|\CFfbi^0\br|\mi 2\br)\1t.
\leqno(3.4.3)$$
\par\nin Here $-2$ after $\bl|\CFfbi^0\br|$ comes from the last condition in the definition of $m_{\si}$.
Note that the $f_{(i)}$ ($i\ins I_f$) are {\it convenient\1} by the hypothesis on $\Gp(f)$.
\sk
We then see that (\hl{3.4.2}{3.4.2}) can be expressed in the form of the right-hand side of (\hl{1.3.1}{1.3.1}) using the calculation in \hl{1.4}{1.4} (for instance, $\be_{i,j}$ is given as in (\hl{1.4.3}{1.4.3})).
Note that $\gat_{\si^{(i)}}-3$ for $\si^{(i)}$ coincides with $\bl|\CFfbi^0\br|\mi 2$, since these are respectively equal to
$$\bl|\CFflinci^2\br|-1=\bl|\CFfbi^1\br|\mi1.$$
\par\nin (Here $\bl|\CFfbi^0\br|$ can be bigger than the number of interior 0-dimensional faces of $\Gp(f_{(i)})$ if there is $\si'\ins\CFfb^1$ such that its image by $\pi_i$ is contained in the boundary of $\Gp(f_{(i)})$, but does not coincide with any 1-dimensional face of $\Gp(f_{(i)})$, that is, if a fusion of compact faces under a projection occurs.)
The equality (\hl{9}{9}) then follows from Theorem\,\,\hl{T1}{1} and the Yomdin-Steenbrink formula for spectrum (\hl{1.3.1}{1.3.1}). This finishes the proof of Theorem\,\,\hl{T3}{3}.
\par\htt{R3.4}{}\msn
{\bf Remark\,\,3.4.} The Yomdin-Steenbrink formula for spectrum (\hl{1.3.1}{1.3.1}) {\it cannot\1} be extended to {\it spectral pairs\1} if there is a {\it fusion of compact faces under a projection\1} as in Remark\,\,\hl{R1.4a}{1.4a}, see also \cite{KL1}, \cite{KL2}.
Indeed, set
$$f_u\eq u\1x^2y^2\pl x^4z\pl y^4z\pl x^5\pl y^6\q(u\ins\De_{\ep}),$$
\par\nin for $0\slt\ep\,{\ll}\,1$. It has non-isolated singularities along the $z$-axis, and its general hyperplane section $h$ has an ordinary quadruple point (that is, $h$ is a $\mu$-constant deformation of $v^4\pl w^4$). Hence its spectral numbers are given by
$$(i{+}j)/4\,\,\,\,(i,j\ins[1,3]),\,\,\,\,\h{that is,}\,\,\,\,\,\Sp_h(t)\eq(t^{1/4}\pl t^{1/2}\pl t^{3/4})^2,$$
\par\nin in particular $\al_{1,1}\eq\tfrac{1}{2}$. Note that the non-integral spectral numbers of $h$ have the {\it middle\1} weight $1\,({=}\,(n{-}1){-}1$), since the Milnor monodromy of $h$ is {\it semisimple.}
\sk
Consider $f_u{+}\1\ell^r$ ($r\,{\gg}\,0$) with $\ell$ a sufficiently general linear function. The Newton polytope of $f_u{+}\1\ell^r$ has an interior 1-dimensional compact face $\sigma$ joining $(2,2,0)$ and $(0,0,r)$ if $u\nes0$, and there is no $r$-depending interior\ 1-dimensional compact face of $\Gp(f_u{+}\1\ell^r)$ if $u\eq0$.
$$\h{$\setlength{\unitlength}{.4cm}
\begin{picture}(15,13.3)
\put(5,5){\vector(-1,-1){4.3}}
\put(5,5){\vector(1,0){8}}
\put(5,5){\vector(0,1){8}}
\qbezier(12,5)(10.5,5.5)(9,6)
\qbezier(1.4,1.4)(2,2.5)(2.6,3.6)
\qbezier(2.6,3.6)(3.8,7.8)(5,12)
\qbezier(9,6)(7,9)(5,12)
\qbezier(5.8,3.8)(5.4,7.9)(5,12)
\qbezier(5.8,3.8)(3.6,2.6)(1.4,1.4)
\qbezier(5.8,3.8)(4.2,3.7)(2.6,3.6)
\qbezier(5.8,3.8)(7.4,4.9)(9,6)
\qbezier(5.8,3.8)(8.9,4.4)(12,5)
\put(5.8,3){$\scriptstyle(2,2,0)$}
\put(10.8,4.1){$\scriptstyle(0,6,0)$}
\put(2,1.2){$\scriptstyle(5,0,0)$}
\put(9.3,6){$\scriptstyle(0,4,1)$}
\put(0.3,3.6){$\scriptstyle(4,0,1)$}
\put(5.3,12){$\scriptstyle(0,0,r)$}
\put(5,0){$\scriptstyle(u\nes 0)$}
\end{picture}$
$\setlength{\unitlength}{.4cm}
\begin{picture}(15,13.3)
\put(5,5){\vector(-1,-1){4.3}}
\put(5,5){\vector(1,0){8}}
\put(5,5){\vector(0,1){8}}
\qbezier(12,5)(10.5,5.5)(9,6)
\qbezier(1.4,1.4)(2,2.5)(2.6,3.6)
\qbezier(12,5)(6.7,3.2)(1.4,1.4)
\qbezier(9,6)(5.8,4.8)(2.6,3.6)
\qbezier(2.6,3.6)(3.8,7.8)(5,12)
\qbezier(2,2)(5.5,4)(9,6)
\qbezier(9,6)(7,9)(5,12)
\put(10.8,3.9){$\scriptstyle(0,6,0)$}
\put(2,1){$\scriptstyle(5,0,0)$}
\put(9.3,5.9){$\scriptstyle(0,4,1)$}
\put(0.3,3.6){$\scriptstyle(4,0,1)$}
\put(5.3,12){$\scriptstyle(0,0,r)$}
\put(5,0){$\scriptstyle(u\eq 0)$}
\end{picture}$}$$
\par\nin \sk
Assume $r$ is a {\it prime\1} at lest 7. Put $r'\defs\tfrac{r{-}1}{2}$. We see that the parallelotope spanned by $(2,2,a)$ and $(0,0,r)$ contains the lattice points $(1,1,k)$ for $k\ins[1,r]$ ($a=0,1$), and the polynomial $f_u{+}\1\ell^r$ has the following {\it non-integral $r$-depending\1} spectral numbers:
\htt{3.4.4}{}
$$\bl\{\tfrac{1}{2}\pl\tfrac{k}{r}\br\}_{k\in K}\,\,(u\nes0),\q\bl\{\tfrac{1}{2}\pl\tfrac{k-1/2}{r}\br\}_{k\in K'}\,\,(u\eq0),
\leqno(3.4.4)$$
\par\nin where $K\defs\{1,\dots,r{-}1\}$, $K'\defs\{1,\dots,r\}\stm\{r'{+}1\}$. Indeed, the planes containing $(2,2,0)$, $(4,0,1)$, $(0,0,r)$ and $(4,0,1)$, $(0,4,1)$, $(0,0,r)$ are defined respectively by
$$ (r{-}1)x{+}(r{+}1)y{+}4z\eq4r,\q\q(r{-}1)x{+}(r{-}1)y{+}4z\eq4r.$$
\par\nin (These can be used to compute the monodromy zeta functions, see \cite{Va1} or \hl{1.5}{1.5}.) Here the {\it non-integral\1} spectral numbers $\tfrac{1}{2}\pl\tfrac{k}{r}$ ($k\ins K$) for $u\nes0$ must have weight different from the middle weight $n{-}1\eq2$, which is 1 in this paper (but may be 3 in another normalization, see \hl{3.1}{3.1}), since the corresponding lattice points are contained in the cone of an {\it interior\1} $1$-dimensional compact face. However, this cannot be compatible with the Yomdin-Steenbrink formula for spectral pairs as long as the spectral number $\al_{1,1}\eq\tfrac{1}{2}$ is endowed with the {\it middle weight\1} 1. If we change the formula itself for this case, it does not work for the {\it non-integral\1} spectral numbers of $f_0{+}\1\ell^r$ or those of $f_u{+}\1\ell^r$ ($u\ne0$) corresponding to lattice points contained in the {\it interiors\1} of the cones of 2-dimensional compact faces. Here the problem is that the image of the cone of the {\it interior\1} 1-dimensional compact face of $f_u{+}\1\ell^r$ by the projection $\pi_3$ is {\it not\1} contained in the cone of a vertex of $\Gp(f_{(3)})$ because of a fusion of compact faces under the projection as is explained below. Note that $\al_{1,1}\eq\be_{1,1}\eq\tfrac{1}{2}$ for $u\eq0$, but this does not hold for $u\nes0$, since $\{f_u|_{z=c}\}$ {\it is not a $\mu$-constant family at\1} $c\eq\pm u/2$ (with $u\nes0$ fixed). From the formula, we see that $\be_{1,1}\eq0$ for $u\nes0$, since the third coordinate of $(2,2,0)$ is 0, see also Question\,\,\hl{Q1.4}{1.4} and Proposition\,\,\hl{P1.4}{1.4}.
\sk
The problem is that the Newton polytope of a general hyperplane section of $f_u$ has a unique 1-dimensional compact face $\tau$ joining $(4,0)$ and $(0,4)$, but this is the {\it union\1} of the images of two 1-dimensional compact faces $\tau'$, $\tau''$ of the Newton polytope of $f_u$ ($u\nes0$) under the projection $\pi_3\,{:}\,\R^3\tos\R^2$ forgetting the third coordinate with $\tau'$, $\tau''$ the convex hulls of $\{(2,2,0),(4,0,1)\}$ and $\{(2,2,0),(0,4,1)\}$ respectively, that is, a {\it fusion of compact faces\1} under a projection occurs if $u\nes0$, see also Remark\,\,\hl{R1.4a}{1.4a}.
\bs\bs\htt{A}{}
\vbox{\centerline{\bf Appendix: Symmetry of combinatorial polynomials}
\ms
\centerline{by Morihiko Saito}
\ms
\centerline{(This was originally contained in \cite{des}.)}}
\bsn
In the {\it simplicial convenient\1} case, we have the {\it symmetry of combinatorial polynomials\1} $r_{\!\si}(t)$ as in (\hl{6}{6}) (which is deduced from the one for the spectrum $\Sp_f(t)$ in (\hl{1.1.5}{1.1.5})) as follows.
\par\htt{TA}{}\msn
{\bf Theorem\,\,A.} {\it If $f$ is Newton non-degenerate with simplicial convenient Newton polytope, then we have the symmetry of combinatorial polynomials}
\htt{A.1}{}
$$r_{\!\si}(t)\eq r_{\!\si}(t^{-1})\1t^{n-d(\si)}\q(\emptyset\less\forall\,\si\less\Gf).
\leqno{\rm(A.1)}$$
\par\nin \msn
{\it Proof.} Let $v_j$ ($j\ins[1,r]$) be the vertices of $\Gf$. For $\si\less\Gf$, let $J(\si)\sst[1,r]$ be the subset such that $j\ins J(\si)$ if and only if $v_j$ is a vertex of the simplex $\si$. There is a positive integer $b$ such that the following inclusions hold for any $\si\less\Gf$\,:
\htt{A.2}{}
$$\msum_{j\in J(\si)}\,\N\1v_j\,\subset\,C(\si)\cap\N^n\,\subset\,\msum_{j\in J(\si)}\,\N\1b^{-1}v_j,
\leqno{\rm(A.2)}$$
\par\nin where $C(\si)\sst\R^n$ is the cone of $\si$. The first inclusion is trivial, since $v_j\ins\N^n\,({=}\,\Z_{\ges 0}^n)$. For the second inclusion, let $V(\si)\sst\R^n$ be the vector subspace spanned by $\si$. We have an inclusion of lattices
$$V(\si)\cap\Z^n\sst\msum_{j\in J(\si)}\,\Z\1b^{-1}v_j,$$
\par\nin for some positive integer $b$, which can be taken independently of $\si$ (since there are finitely many $\si$).
The second inclusion follows from this, since
$$C(\si)\eq\msum_{j\in J(\si)}\,\R_{\ges 0}\1v_j.$$
\par\nin \sk
Let $p_j$ ($j\ins[1,r]$) be mutually distinct primes with $p_j\gg nb$. We assume $p_1\eq\min\{p_j\}$. Set
$$g\eq\msum_{j=1}^r\,x^{p_jv_j}.$$
\par\nin Since $f$ is {\it simplicial,} there is a very small positive number $\ep_f$ such that the combinatorics of the compact faces of $\Gp(g)$ is the same as that of $\Gp(f)$ if the following condition is satisfied\,:
\htt{A.3}{}
$$0\slt(p_j\mi p_1)/p_1\slt\ep_f\q(\forall\,j\ins[2,r]).
\leqno{\rm(A.3)}$$
\par\nin By Lemma\,\,\hl{LA}{A} below, we may assume that condition (\hl{A.3}{A.3}) is satisfied (replacing the $p_j$ if necessary).
Proposition\,\,\hl{PA}{A} below then implies the assertion (\hl{A.1}{A.1}) using (\hl{5}{5}) and the self-dualities\,:
\htt{A.4}{}
$$\Sp_g(t)\eq\Sp_g(t^{-1})\1t^n,\q q_{\si}(t)\eq q_{\si}(t^{-1})\1t^{d(\si)},
\leqno{\rm(A.4)}$$
\par\nin see (\hl{1.1.5}{1.1.5}) for the first equality. This finishes the proof of Theorem\,\,\hl{TA}{A} (admitting Lemma\,\,\hl{LA}{A} and Proposition\,\,\hl{PA}{A} below).
\par\htt{PA}{}\msn
{\bf Proposition\,\,A.} {\it With the above notation and assumption, the $q_{\si}(t)$ for $\si\less\Ga_{\!g_a}$ $($including the case $\si\eq\emptyset\1)$ are linearly independent over $\C(\!(t)\!)\,({=}\,\C[[t]][t^{-1}])$ in $\C[[t^{1/e}]][t^{-1}]$, where $e$ is a sufficiently divisible positive integer.}
\msn
{\it Proof.} For $\nu\eq\msum_{j\ins J(\si)}\,(k_j/b)v_j\ins C(\si)^{\circ}\cap\N^n$ with $k_j\ins\Z_{>0}$, we have
\htt{A.5}{}
$$v'_g(x^{\nu})\eq\msum_{j\in J(\si)}\,k_j/p_jb,
\leqno{\rm(A.5)}$$
\par\nin since $v'_g(x^{p_jv_j})\eq1$ (and $\si$ is a simplex). Write
$$v'_g(x^{\nu})\eq A_{\nu}/B_{\nu},$$
\par\nin with $A_{\nu},B_{\nu}$ mutually prime positive integers. Then we have
\htt{A.6}{}
$$\aligned&B_{\nu}\eq b_{\nu}\1\mprod_{j\in J(\si)}\,p_j\,\,\,\,\,\h{with}\,\,\,\,\,b_{\nu}\,|\,b\,\,\,\,\,\h{(that is,}\,\,\,\,b/b_{\nu}\ins\Z),\\&\h{if}\,\,\,\,\,0\slt k_j\less nb\,\,\,\,(\forall\,j\ins J(\si)).\endaligned
\leqno{\rm(A.6)}$$
\par\nin Indeed, it is easy to see that $B_{\nu}\,|\,b\1\bl(\mprod_j\,p_j\br)$, and for any $i\ins J(\si)$, we have
$$b\1\bl(\mprod_{j=1}^r\,p_j\br)\bl(\msum_{j\in J(\si),\,j\ne i}\,k_j/p_jb\br)\in\Z\1p_i.$$
\par\nin This implies that, if $p_i\nmid B_{\nu}$ (that is, $B_{\nu}/p_i\notin\Z$) so that
$$b\1\bl(\mprod_{j=1}^r\,p_j\br)\,A_{\nu}/B_{\nu}\in\Z\1p_i,$$
\par\nin then we get that $\bl(\mprod_{j\ne i}\,p_j\br)\,k_i\ins\Z\1p_i$, hence $p_i\,|\,k_i$. However, this contradicts the assumption that $0\slt k_j\less nb$, since $nb\ll p_i$. So (\hl{A.6}{A.6}) follows.
\sk
We denote by ${\rm Supp}\,q_{\si}(t)\sst\Q$ the support of $q_{\si}(t)$, that is, the set of $\al\ins\Q$ with $q_{\si,\al}\ne 0$, where $q_{\si}(t)\eq\msum_{\al}\,q_{\si,\al}\1 t^{\al}$. By (\hl{A.6}{A.6}) we get that
\htt{A.7}{}
$${\rm Supp}\,q_{\si}(t)\cap{\rm Supp}\,q_{\si'}(t)\cap[0,n/p_1)=\emptyset\q\h{if}\q\si\ne\si'.
\leqno{\rm(A.7)}$$
\par\nin Indeed, we have
$${\rm Supp}\,q_{\si}(t)\eq\{v'_g(x^{\nu})\mid\nu\ins E_{\si}^{\circ}\cap\N^n\},$$
\par\nin where $E_{\si}^{\circ}$ is the interior of the parallelotope $E_{\si}$ spanned by the vertices of the simplex $\si$ if $\si\ne\emptyset$ (and ${\rm Supp}\,q_{\si}(t)\eq\{0\}$ for $\si\eq\emptyset$). Moreover, the condition $v'_g(x^{\nu})\slt n/p_1$ implies that $k_j\less nb$ with $k_j$ as in (\hl{A.5}{A.5}) (using the second inclusion of (\hl{A.2}{A.2})). So (\hl{A.7}{A.7}) follows from (\hl{A.6}{A.6}).
\sk
By the first inclusion of (\hl{A.2}{A.2}), we also see that
\htt{A.8}{}
$${\rm Supp}\,q_{\si}(t)\cap[0,n/p_1)\ne\emptyset\q\q(\emptyset\less\forall\,\si\less\Gg).
\leqno{\rm(A.8)}$$
\par\nin Indeed, setting $v_{\si}\defs\sum_{j\in J(\si)}v_j\ins C(\si)^{\circ}\cap\Z^n$, we have $v'_f(v_{\si})\eq\msum_{j\in J(\si)}\,1/p_j\slt n/p_1$, since $p_1\slt p_j$ ($\forall\,j\ins[2,r]$).
\sk
The linear independence now follows form (\hl{A.7}{A.7}--\hl{A.8}{8}). Indeed, assume
$$\msum_{\emptyset\les\si\les\Gg}\,c_{\si}(t)\1q_{\si}(t)\eq0\q\h{with}\q c_{\si}(t)\ins\C(\!(t)\!).$$
\par\nin Here we may assume $c_{\si}(t)\ins\C[[t]]$ ($\forall\,\si$) and $c_{\si}(0)\ne0$ ($\exists\,\si$), since we can divide the $c_{\si}(t)$ by a power of $t$ (independently of $\si$). Then this equality contradicts (\hl{A.7}{A.7}--\hl{A.8}{8}), and Proposition\,\,\hl{PA}{A} follows.
\par\htt{LA}{}\msn
{\bf Lemma\,\,A.} {\it For any positive number $\ep$ and positive integers $m,r$, there are mutually distinct prime numbers $\,p_j\gess m\,$ $(j\ins[1,r])$ such that $0\slt(p_j\mi p_1)/p_1\slt\ep$ for $j\ins[2,r]$.}
\msn
{\it Proof.} If the assertion does not hold, the number of primes contained in the interval
$$[m(1{+}\ep)^i,m(1{+}\ep)^{i+1})$$
\par\nin must be at most $r{-}1$ for any $i\ins\N$. This implies that
\htt{A.9}{}
$$\aligned\pi(x)&\less\pi(m)\pl(r{-}1)\bl\lceil(\log x\mi\log m)/\log(1{+}\ep)\br\rceil\\&\slt C_1\pl C_2\log x\q\q\q(\forall\,x\gess m).\endaligned
\leqno{\rm(A.9)}$$
\par\nin Here $C_1,C_2\ins\R_{>0}$, $\lceil\al\rceil\defs\min\{k\ins\Z\mid\al\less k\}$, and $\pi(x)$ denotes the number of primes which are at most $x\ins\R_{>0}$.
(Indeed, if $x/m\eq(1{+}\ep)^u$, then $u\eq(\log x\mi\log m)/\log(1{+}\ep)$.)
However, this inequality contradicts the {\it prime number theorem\1} claiming that
\htt{A.10}{}
$$\pi(x)\sim x/\log x\q\bl(\h{that is,}\,\,\,\lim_{x\to+\infty}\pi(x)\log x/x\eq1\br).
\leqno{\rm(A.10)}$$
\par\nin (Put $y\eq\log x$ so that $x/\log x\eq e^y/y$, and consider $\rho(y)\defs\pi(e^y)$.) Thus Lemma\,\,\hl{LA}{A} is proved.
\par\htt{RA}{}\msn
{\bf Remark\,\,A.} It is not necessarily clear whether the coefficients of $\SpDef_f(t)$ are non-negative. More precisely, writing the combinatorial polynomials $r_{\!\si}(t)$ as $\msum_k\,r_{\si,k}\1t^k$ with $r_{\si,k}\ins\Z$, it is quite unclear whether the coefficients $r_{\si,k}$ are weakly increasing for $k\less(n\mi d(\si))/2$, and weakly decreasing for $k\gess(n\mi d(\si))/2$, although this could be expected very much by the monodromical property of the weight filtration, see (\hl{1.1.8}{1.1.8}). This is closely related to an observation that $q_{\si}(t)\mi q_{\tau}(t)\1t^i$ have non-negative coefficients when $d_{\si}\mi d_{\tau}\eq 2i\sgt0$ as far as calculated. In the case $n\eq4$, it does not seem clear whether $\ep_{\emptyset}\gess 0$.

\end{document}